# Structured $H_\infty$-control of infinite dimensional systems

P. Apkarian[1]*, D. Noll[2]

**Abstract**

We develop a novel frequency-based $H_\infty$-control method for a large class of infinite-dimensional Linear-Time-Invariant systems in transfer function form. Major benefits of our approach is that reduction or identification techniques are not needed thereby avoiding possible distortions. It can exploit either transfer function models or input/output frequency response data when available. It computes simple practical controllers of any size and structure.

We use a non-smooth trust-region method to compute arbitrarily structured locally optimal $H_\infty$-controllers for a frequency sampled approximation of the underlying infinite-dimensional $H_\infty$-problem in such a way that exponential stability in closed-loop is guaranteed, and the optimal value of the approximation captures the value of the infinite-dimensional problem within a prior tolerance level. We demonstrate the versatility and practicality of our method on a variety of infinite-dimensional $H_\infty$-synthesis problems, including distributed and boundary control of PDEs, control of dead time and delay systems, and using a rich testing set.

[1] Control System Department, ONERA, 2, av. Ed. Belin, 31055, Toulouse, France
[2] Institut de Mathématiques de Toulouse, 118 route de Narbonne F-31062, Toulouse, France
*Corresponding author: P. Apkarian & D. Noll

## Contents

## Introduction

In this work we use a frequency-based $H_\infty$-method to control infinite-dimensional LTI-systems $G(s)$. After embedding $G(s)$ as usual in a plant $P(s)$ and setting up performance and robustness channels $T_{wz}(P,K)$, we address the infinite-dimensional $H_\infty$-optimization problem

$$
\begin{aligned}
\text{minimize} \quad & \max_{\omega \in [0,\infty]} \overline{\sigma}\left(T_{zw}(P(j\omega), K(j\omega))\right) \\
\text{subject to} \quad & K \text{ stabilizes } G \text{ exponentially} \\
& K \in \mathscr{K}
\end{aligned}
\tag{1}
$$

where optimization is over a class $\mathscr{K}$ of structured finite rank control laws. Our strategy is to choose frequency samples $G(j\omega_\nu)$ of $G(s)$ in such a way that the solution $K^* \in \mathscr{K}$ of the approximate $H_\infty$ program

$$
\begin{aligned}
\text{minimize} \quad & \max_{\nu=1,\dots,N} \overline{\sigma}\left(T_{zw}(P(j\omega_\nu), K(j\omega_\nu))\right) \\
\text{subject to} \quad & K \text{ stabilizes } G \text{ exponentially} \\
& K \in \mathscr{K}
\end{aligned}
\tag{2}
$$

guarantees closed-loop stability for $G(s)$, and assures that the value of (2) differs only by a fixed tolerance $\vartheta$ from the true value of (1). Sampling in the frequency domain becomes necessary since the objective of (1) is semi-infinite, non-smooth, and non-convex, and not directly amenable to efficient computation.

The difficulty in program (1) is further aggravated by the fact that controllers $K \in \mathscr{K}$ have to be *structured* in the sense of [1]. Structured controllers or control architectures are preferred by practitioners and include classics like PIDs, lead-lag and notch filters, polynomial matrix fractions, reduced fixed-order controllers, but also observer-based controllers, distributed control architectures interconnecting other structured elements, decentralized control, and much else. A general way to model structure uses state-space in the form

$$K(s): \begin{cases} \dot{x}_K &=& A_K(\mathbf{x})x_K + B_K(\mathbf{x})y \\ u &=& C_K(\mathbf{x})x_K + D_K(\mathbf{x})y. \end{cases} \quad (3)$$

where $A_K(.)$, $B_K(.)$, ... are smooth matrix-valued functions of a tunable parameter vector $\mathbf{x} \in \mathbb{R}^n$, but our method applies also to infinite-dimensional controller structures $K(\mathbf{x})$ as long as they are parametrized by a finite-dimensional vector $\mathbf{x} \in \mathbb{R}^n$ of tunable parameters. With this restriction programs (1), (2) fall within the class of semi-infinite optimization problems [2].

For systems given in state-space we obtain the transfer function $G(s)$ directly from the infinite dimensional system. We then discretize $G(j\omega_\nu)$ on the low-dimensional level of the input-output map, where sampling is best adapted to the truly relevant dynamics of the system. In consequence, versions of (2) with essentially no loss over (1) typically require only a very moderate number $N$ of samples, rarely exceeding a couple of hundred nodes, so that (2) is solved in seconds to minutes. Pre-computing samples $G(j\omega_\nu)$ may turn out more time-consuming, but since we perform it offline, it does not impede the optimization or the plant-modeling phase. Our method is also suited for systems provided from start in frequency sampled form (2), or for systems given directly by their transfer function.

In order to justify our approach theoretically, we have to clarify the following issues:

(a) How to sample the return difference $\det(I + GK)$ so that exponential stability in closed loop is guaranteed.

(b) How to sample the transfer function $G(s)$ so that the approximate value of (2) is within a fixed tolerance $\vartheta$ of the true value of (1).

(c) How to solve the non-smooth optimization program (2) algorithmically.

To address the stability issue (a) we implement an infinite-dimensional Nyquist test, which is effective as soon as stability of the closed-loop systems is *spectrum-determined*. This applies for instance to delay and dead-time systems, to boundary and distributed control for parabolic PDEs, or to control of hyperbolic PDEs of one space dimension. In contrast, control of hyperbolic PDEs of several space dimensions requires a case-by-case analysis.

Analysis of the stability issue (a) reveals the surprising fact that most of the time only a very limited number of frequency samples $G(j\omega_\nu)$ are needed to obtain the correct winding number. On the other hand, the sampling grid for stability $\Omega_{\text{nyq}}$ depends on the candidate controllers $K \in \mathscr{K}$, and therefore needs updating during optimization.

A second important aspect of question (a) is how stability can be built into a mathematical programming constraint in order to maintain it during optimization (2). This cannot be based on the Nyquist test, which is discrete in nature, and we propose a stability barrier function based on the modulus margin of the closed-loop system.

Appropriate frequency sampling to assure performance (b) benefits from the fact [3] that for a fixed control law $K \in \mathscr{K}$ the frequency response, even when exhibiting sharp primary and secondary peaks, is twice continuously differentiable as a function of frequency in the neighborhood of those peaks. This improves the order of the approximation and leads to an efficient sampling method. Non-smoothness typically occurs at anti-resonances, but as those are irrelevant for a good approximation of the frequency maximum in (1), the number of nodes needed for a good approximations is moderate and rarely exceeds a couple of hundred.

Optimization (c) is based on the non-smooth trust-region method of [4], which was already successfully applied to computation of the structured distance to instability in [5]. For the solution of program (2) specific features of the method are exploited to gain speed, on which we comment in section 4. A technical difficulty arises from the fact that the sampling grid for performance $\Omega_{\text{opt}}$, unlike for stability, cannot be adapted to the candidate controllers $K \in \mathscr{K}$ during optimization, as this would change program (2). Posterior verification of the optimal controller $K^* \in \mathscr{K}$ obtained on the current $\Omega_{\text{opt}}$ is therefore necessary, and may require occasional restarts on a refined optimization grid. The overall procedure including these re-starts given in algorithm 3 is still speedy and converges within seconds to minutes.

A side aspect of our approach is that it avoids the use of system reduction and identification techniques, and allows us to stay as close as possible to the infinite-dimensional program (1). Discretization if any is confined to the level of the transfer function.

Even though our primary interest here lies in situations where the transfer function $G(s)$ is available analytically, or numerically at arbitrary frequencies, program (2) also contributes novel aspects in cases where from start only a frequency sampled version $G(j\omega_\nu)$ is available for synthesis, with no recourse to further missing values $G(j\omega)$. Our technique may then still be applied to reduce program (2) from the original fine sampling $\Omega_{\text{fine}}$ to a manageable size $\Omega_{\text{opt}}$ for optimization, with stability and performance certificates then valid under the proviso that the information stored in the initial finest available sample $\Omega_{\text{fine}}$ is sufficiently rich.

There is a vide literature on controller design based on frequency-domain data, and we just cite a few. Pioneering work is the semi-infinite programming technique proposed by Polak [6], and the Quantitative Feedback Theory (QFT) of [7] is in this class. More recently, various optimization-based techniques have been studied. Linear programming or convex optimization is proposed for specific controller structures in [8, 9, 10]. A more general convex-concave procedure (CCP) is used in [11] to design PIDs, and in [12] is extended to linearly parameterized MIMO PIDs. In the same vein, the arXiv paper [13] applies CCP to synthesize MIMO fraction-of-polynomial controllers. These specific controller structures allow design specifications in the form of convex differences. Linearizing concave terms then yields LMI subproblems, which are solved sequentially to determine locally optimal controllers. A general analysis of CCP together with variations and extensions is discussed in [14].

Nyquist stability for infinite dimensional systems has a long history and is discussed in [15], an axiomatic approach being [16]. In [17] extensions to trace class operators are proposed. The link between input-output and exponential stability is discussed in [18].

The structure of the paper is as follows. After some preparations in section 1, we discusses theoretical and practical aspects of the Nyquist test in section 2, and its use to enforce closed-loop stability in (2). Grid selection for the Nyquist test is presented in section 3. Section 4 presents our optimization method for (2), grid selection for optimization $\Omega_{\mathrm{opt}}$ being discussed in section 5. Sections 7, 8 discuss control of crystallization and dead-time processes. A numerical evaluation of our method using the test bench [19], along with several PDE studies, is presented in Section 10.

## Notation

Notions from classical control theory are covered by [20], basics on infinite-dimensional systems are found in [21, 22], more details on well-posed systems will be given in the next section. The index of a curve $\gamma$ around a point $x$ is $\mathrm{ind}(\gamma, x)$, see e.g. [23, p. 139]. For a complex valued function $f$, we use $\Delta_{[\omega_1, \omega_2]} \arg f$ to denote the variation of argument of $f$ along the segment $[j\omega_1, j\omega_2]$ of $j\mathbb{R}$. For concepts in non-differential optimization we refer to [24, 4].

## 1. Well-posed transfer functions

We consider well-posed transfer functions $G(s)$, which are generated by well-posed linear systems $\Sigma = (A, B, C, \mathscr{G})$ in the sense of Salomon [25] and Weiss [26], see also Curtain [27]. Here $A$ is the generator of a strongly continuous semigroup on a Hilbert space $X$, $X_1 \subset X$ is $D(A)$ equipped with the graph norm, $X_{-1}$ is the Hilbert space obtained by completing $X$ with respect to the norm $\|x\|_{-1} = \|(\beta I - A)^{-1} x\|$, where $\beta \in \rho(A)$ is fixed, so that $X_1 \subset X \subset X_{-1}$, $B \in L(U, X_{-1})$ is the control operator, $C \in L(X_1, Y)$ is the observation operator, $\mathscr{G} : L^2_\sigma([0, \infty), U) \to L^2_\sigma([0, \infty), Y)$ for some $\sigma \in \mathbb{R}$ is the input-output map, a bounded causal time-invariant linear operator. The transfer function $G(s) \in \mathbf{H}^\infty_\sigma$ is defined on $\mathbb{C}^+_\sigma = \{s \in \mathbb{C} : \mathrm{Re}(s) > \sigma\}$ with values in $L(U, Y)$, which satisfies $\widehat{y} = G\widehat{u}$ whenever $y = \mathscr{G}u$, $u \in L^2_\sigma([0, \infty), U)$. We assume throughout that $G$ is matrix-valued, which means that input and output spaces $U \simeq \mathbb{R}^p$ and $Y \simeq \mathbb{R}^m$ are finite-dimensional. This hypothesis is necessary to assure that the computed control laws are implementable. The transfer function is proper if $\overline{\sigma}(G(s)) \le M$ for some $M > 0$, some $\rho > 0$, and all $s \in \{s \in \mathbb{C}^+ : |s| \ge \rho\}$.

The well-posed system $\Sigma$ is regular if the limit of $G(s)$ as $s \to \infty$ along the positive real axis exists. In that case the direct transmission $D \in L(U, Y)$ is well-defined, and according to [26, Theorem 1.1] the transfer function $G(s)$ may now be represented as $G(s) = C_L(sI - A)^{-1}B + D$, where $C_L$ is a suitable extension of the operator $C$ obtained by applying the Cesàro summability method to the output operator of $\Sigma$, referred to as the Lebesgue extension of $C$ in [26]. One can also use the $\Lambda$-extension $C_\Lambda$, which uses the Abel summability method instead and satisfies $X_1 \subset D(C_L) \subset D(C_\Lambda) \subset X$, extending $C$ even further. The notion of regularity is convenient in so far as the pointwise representation of the transfer function is now almost identical with the classical case with bounded $B, C$, but otherwise regularity is not essential for the present work.

Static output feedback $T(G, K)$ is defined as follows. An operator $K \in L(Y, U)$ is an admissible static output feedback for the well-posed system $\Sigma$ if $I + \mathscr{G}K$ is invertible in the space $TIC_\sigma(U)$ of causal time-invariant operators $L^2_\sigma(\mathbb{R}, U) \to L^2_\sigma(\mathbb{R}, Y)$ for some $\sigma \in \mathbb{R}$. Equivalently this means that $I + G(s)K$ is invertible on $\mathrm{Re}(s) > \sigma$, and its inverse $T(G, K)$ is a well-posed transfer function, see [26, 28, 27].

If $G$ is regular, then the closed-loop transfer function $T(G, K) = (I + G(s))^{-1}$ is also regular. This is a consequence of $\dim(U) < \infty$ and $\dim(Y) < \infty$, see [28, Prop. 4.6], and also [27, p. 216].

Dynamic feedback is introduced as follows. We consider an infinite-dimensional controller $K$ represented in just the same way by a well-posed system with generator $A_K$ on a Hilbert space $X^K$, control operator $B_K \in L(Y, X^K_{-1})$ with input space $Y \simeq \mathbb{R}^m$, observation operator $C_K \in L(X^K_1, U)$ with output space $U \simeq \mathbb{R}^p$, and transfer function $K(s)$. Now we define the lower feedback connection $T(G, K)$ by forming the cross product $G \times K$ (also known as the parallel connection) of system and controller, saying that $K$ is an admissible dynamic lower output feedback for $G$ if the static operator $J := \begin{bmatrix} 0 & I \\ -I & 0 \end{bmatrix}$ is an admissible static output feedback for $G \times K$ in the sense introduced above, see e.g. [29, Theorem 3.4]. In other words, $T(G, K) := T(G \times K, J)$, these definitions being consistent when $K$ is static.

Since the cross product has transfer function $\mathrm{diag}(G, K)$, we find that admissibility of the dynamic feedback requires that $F(s) := I \times I + \mathrm{diag}(G(s), K(s))J$ be invertible and its inverse $T(s) = F(s)^{-1}$ be a well-posed transfer function. Writing more explicitly

$$F(s) = \begin{bmatrix} I & G(s) \\ -K(s) & I \end{bmatrix} \tag{4}$$

we find that its inverse is

$$T = \begin{bmatrix} I - K(I + GK)^{-1}G & -K(I + GK)^{-1} \\ (I + GK)^{-1}G & (I + GK)^{-1} \end{bmatrix}. \tag{5}$$

If $G, K$ are regular, then so is $G \times K$, and it follows from the above that the closed-loop system $T(G \times K, J)$ is automatically regular. In the regular case state-space representations of the closed loop then resemble those known in the finite-dimensional case, and we refer to [27, 26, 30] for more details and explicit formulas.

The well-posed system $\Sigma$ is internally stable if $A$ generates an exponentially stable semigroup. The system is externally stable if its transfer function $G(s)$ belongs to the space $\mathbf{H}_\sigma^\infty$, with $\sigma = 0$, which we abbreviate by $\mathbf{H}^\infty$. In closed loop, external stability is therefore expressed as $T \in \mathbf{H}^\infty$.

Following Morris [29], $(A,B)$ is exponentially stabilizable if there exists an observation operator $K$ such that the system $(A,B,K,\mathscr{G}_K)$ is well-posed, admits $-I$ as an admissible static feedback operator, and the closed loop is internally stable. Exponential detectability is defined analogously.

**Lemma 1.** (Morris [29, Thm. 5.2], see also Rebarber [31], Curtain [32]). *A well-posed system is exponentially stable if and only if it is exponentially stabilizable, exponentially detectable, and externally stable.*

## 2. Winding number, Nyquist stability

Given a class $\mathscr{K}$ of admissible dynamic controllers for $G$ and some candidate $K \in \mathscr{K}$, we define $F(s)$ as in (4) and let $f(s) = \det F(s) = \det(I + G(s)K(s))$. When $F(s)$ is meromorphic on a domain containing $\overline{\mathbb{C}}^+$, we define $n_p$ as the number of poles of $F$ in $\mathbb{C}^+$. We need the following hypotheses:

(i) $G$ and $K$ are proper.

(ii) $F$ has no zeros on $j\mathbb{R}$.

(iii) The limit of $f(s) = \det(I + G(s)K(s))$ as $s \to \infty$ on $\overline{\mathbb{C}}^+$ exists and differs from 0.

(iv) The realizations of $G$ and $K \in \mathscr{K}$ are exponentially stabilizable and detectable.

It follows from (i) that $F(s)$ has only finitely many poles on $j\mathbb{R}$. Now let $h$ be a holomorphic function on a domain containing $\overline{\mathbb{C}}^+$ such that $h(s) \neq 0$ on $\mathbb{C}^+$, $\lim_{s \to \infty} h(s) = 1$ on $\overline{\mathbb{C}}^+$, and such that $h$ has a zero of order $p$ at 0 or $\pm j\omega$ precisely when $F(s)$ has a pole of order $p$ at 0 or $\pm j\omega$. (If $F$ has no poles on $j\mathbb{R}$, then $h = 1$.) Since poles of $f$ are also poles of $F$, $h$ removes also all poles of $f$ on $j\mathbb{R}$. We put $\widetilde{f} = fh$ and call $\{\widetilde{f}(j\omega) : \omega \in \mathbb{R} \cup \{\infty\}\}$ the modified infinite Nyquist curve.

**Theorem 1.** *Let conditions* (i) - (iv) *be satisfied. Suppose the modified infinite Nyquist curve* $\{\widetilde{f}(j\omega) : \omega \in \mathbb{R} \cup \{\infty\}\}$ *winds $n_p$ times around the origin in the clockwise sense, i.e.*

$$\frac{1}{2\pi} \int_{-\infty}^{\infty} \frac{\widetilde{f}'(j\omega)}{\widetilde{f}(j\omega)} d\omega = n_p. \tag{6}$$

*Then the closed-loop system is exponentially stable.*

*Proof.* 1) Since $G$ is exponentially stabilizable and $\dim(U) < \infty$, it follows from Staffans [30, Lem. 8.2.9] that $G(s)$ admits a meromorphic extension on a domain containing $\mathrm{Re}(s) >$ $-\alpha$ for some $\alpha > 0$. Since $K$ is exponentially stabilizable and $\dim(Y) < \infty$, the same is true for $K$. Then $F$ and $f$ are also meromorphic on $\mathrm{Re}(s) > -\alpha$. Since $G, K$ are both proper by condition (i), they have only finitely many poles in $\overline{\mathbb{C}}^+$, and hence so has $F$. In particular, the definition of $n_p$ as the number of poles of $F$ in $\mathbb{C}^+$ makes sense.

2) By hypothesis (iii) the limit of $f(s) = \det(I + G(s)K(s))$ as $s \to \infty$ on $\overline{\mathbb{C}}^+$ exists, and since $f$ is meromorphic by part 1), it has only finitely many poles on the right half plane $\overline{\mathbb{C}}^+$. The same is true for $\widetilde{f}$, and since $h$ removes the poles of $F$ on $j\mathbb{R}$, it also removes the poles of $f$ on $j\mathbb{R}$, so that the number of poles of $\widetilde{f}$ on $\overline{\mathbb{C}}^+$ equals the number of poles of $f$ on $\mathbb{C}^+$. Let us call this number $\widetilde{n}_p$. Moreover, by part 1) we may find a domain $\Omega$ containing $\overline{\mathbb{C}}^+$ on which the number of poles $\widetilde{n}_p$ of $\widetilde{f}$ remains the same. It also follows from (iii) that the modified Nyquist curve is closed.

Next, since by (iii) the limit of $f$ at infinity is different from 0, we infer that $f$ has only finitely many zeros on $\mathbb{C}^+$, and we denote their number by $\widetilde{n}_z$. By (ii) $f$ has no zeros on $j\mathbb{R}$, because zeros of $f$ are also zeros of $F$, hence $f$ has $\widetilde{n}_z$ zeros on a domain $\Omega$ containing $\overline{\mathbb{C}}^+$.

By construction $h$ removes the poles of $F$ on $j\mathbb{R}$ and has no zeros on $\mathbb{C}^+$, so by adjusting $\Omega$ if necessary we may assume that $h$ has no further zeros on $\Omega$. Since $F$ has no zeros on $j\mathbb{R}$ by (ii), there cannot be any cancellations of zeros and poles in computing $f$ on $j\mathbb{R}$, so $h$ not only just removes the poles of $f$ on $j\mathbb{R}$, it also does not add any superfluous zeros on $j\mathbb{R}$. Altogether, $\widetilde{f} = fh$ has therefore neither poles nor zeros on $j\mathbb{R}$, which means it has $\widetilde{n}_z$ zeros on the domain $\Omega$ containing $\overline{\mathbb{C}}^+$. This also shows that the modified Nyquist curve is well-defined and does not pass through the origin, and that the function $\widetilde{f} = fh$ is now amenable to the argument principle on the Nyquist curve with regard to the origin.

3) Consider a standard finite Nyquist D-contour $D$, and for $\varepsilon > 0$ let $D_\varepsilon$ be its $\varepsilon$-enlargement into the left half plane. In other words, while in $D$ we cut the circle to a half-circle along $j\mathbb{R}$, $D_\varepsilon$ corresponds to cutting the circle at $-\varepsilon + j\mathbb{R}$. Suppose $D$ is large enough to contain all rhp poles of $F$, all $\widetilde{n}_p$ rhp poles of $f$, and all $\widetilde{n}_z$ rhp zeros of $f$ in its interior. Assure that $\varepsilon$ is small enough so that $D_\varepsilon$ contains none of the stable poles and zeros of $F, f$, which could arise inside the D-contour on $\mathrm{Re}(s) > -\alpha$, and is contained in $\Omega$. Note that $F, f$ may have infinitely many stable poles and zeros on $\mathrm{Re}(s) > -\alpha$, but only finitely many are within the D-circle, so we may adjust $\varepsilon$ to $D$ to avoid them.

Then by the argument principle the index satisfies

$$\mathrm{ind}(\widetilde{f} \circ D_\varepsilon, 0) = \widetilde{n}_z - \widetilde{n}_p.$$

Passing to the limit $\varepsilon \to 0$ for a fixed D-contour gives

$$\mathrm{ind}(\widetilde{f} \circ D, 0) = \widetilde{n}_z - \widetilde{n}_p,$$

as there are neither zeros nor poles in $D_\varepsilon \setminus D$. By (ii) we have $\lim_{s \to \infty} f(s) \neq 0$ on $\overline{\mathbb{C}}^+$, hence $\lim_{s \to \infty} f(s)h(s) \neq 0$ on $\overline{\mathbb{C}}^+$, and then we may pass to the limit $\tilde{f}(D) \to \tilde{f}(j\mathbb{R})$ in the D-contour to obtain

$$\frac{1}{2\pi} \int_\infty^{-\infty} \frac{\tilde{f}'(j\omega)}{\tilde{f}(j\omega)} d\omega = \tilde{n}_z - \tilde{n}_p,$$

as from a certain size of the D-contour onward the term on the right remains the same. Now by (6), the left hand integral equals $-n_p$, so we have shown $\tilde{n}_p - n_p = \tilde{n}_z$.

Since every rhp pole of $f$ is also a rhp pole of $F$, we have $\tilde{n}_p \leq n_p$. The case $\tilde{n}_p < n_p$ is *a priori* possible and indicates a pole zero cancellation between $G$ and $K$ on $\mathbb{C}^+$. But such a pole zero cancellation would now give $\tilde{n}_z < 0$, which is impossible, as $\tilde{n}_z$ is a natural number. We deduce that $\tilde{n}_p = n_p$, and hence $\tilde{n}_z = 0$, i.e., $\tilde{f}(s)$ has no zeros on $\mathbb{C}^+$, and then neither has $f$. But recall that $\tilde{n}_z$ was also the number of zeros of $\tilde{f}$ within the $D_\varepsilon$ contours for sufficiently large radius and sufficiently small $\varepsilon$, hence $\tilde{f}$ has no zeros on any of those $D_\varepsilon$. Altogether, $\tilde{f}$ has no zeros on a domain $\Omega$ containing $\overline{\mathbb{C}}^+$. Since $f$ has no zeros on $j\mathbb{R}$, the same is true for $f$.

4) We argue that the domain $\Omega$ may be covered by suitable subdomains $\Omega' \subset \Omega$ such that on every $\Omega'$ the matrix function $F(s)$ has a coprime factorizations over the space $\mathbf{H}(\Omega')$ of matrix valued functions holomorphic on $\Omega'$. Choose $\Omega'$ e.g. as a disk contained in $\Omega$ with no poles of $F$ on its boundary, and map it conformally into the hight half $z$-plane $\mathbb{C}^+$ using a Möbius transformation $z = \psi^{-1}(s)$. Then $F \circ \psi$ is meromorphic on $\mathbb{C}^+$, and is bounded as $z \to \infty$ on $\mathbb{C}^+$, because the choice of $\Omega'$ assures that $F(s)$ remains bounded as $s = \psi(z)$ approaches the boundary of $\Omega'$. Hence $\tilde{F}(z) := (z+1)^{-1}F(\psi(z)) = O(z^{-1})$ as $z \to \infty$ on $\mathbb{C}^+$. Therefore by Mossaheb [33] we get a coprime factorization of $\tilde{F}(z)$ over $\mathbf{H}(\mathbb{C}^+)$, which in view of $z + 1 \neq 0$ on $\mathbb{C}^+$ yields a coprime factorization of $F \circ \psi$ over $\mathbf{H}(\mathbb{C}^+)$, and hence via $\psi^{-1}$, a coprime factorization of $F(s)$ over $\mathbf{H}(\Omega')$.

5) Now consider one such $\Omega' \subset \Omega$ and its coprime factorization $F = M^{-1}N$ over $\mathbf{H}(\Omega')$. Since $f = \det(M)/\det(M)$ has no zeros on $\Omega'$, we deduce that neither does $\det(N)$ has no zeros on $\Omega'$. Indeed, from the argument of part 3) we saw that $n_p = \tilde{n}_p$, which meant none of the rhp poles in $F$ disappeared when forming the determinant $f$ due to cancellation with a rhp zero in $F$. But that also means that none of the rhp zeros in $F$ disappears in $f$ due to a cancellation with a rhp pole in $F$. Hence $\det(N)$ has no zeros on $\Omega'$. In other words, we have shown that $f = \det(N)/\det(M)$ is also coprime.

Since $N$ is holomorphic on $\Omega'$ and $\det(N) \neq 0$, it is invertible and its inverse is also holomorphic on $\Omega'$. Then $F^{-1} = N^{-1}M$ is holomorphic on $\Omega'$. But $F^{-1} = T$, where $T$ is the closed-loop transfer function (5), so we have proved that $T$ is holomorphic on $\Omega'$. Since the $\Omega'$ cover $\Omega$, we deduce that $T$ is holomorphic on the domain $\Omega$ containing $\overline{\mathbb{C}}^+$.

6) We argue that $T \in \mathbf{H}^\infty(\mathbb{C}^+)$. For that it remains to prove that $T$ is bounded on $j\mathbb{R}$. But this follows from the fact that any of the four closed-loop transfer functions $G_{cl}$ occurring in $T$ in (5) is proper, i.e. satisfies $\overline{\sigma}(G_{cl}(s)) \leq M$ for some $M > 0$, $\rho > 0$, and all $s$ in $\{s \in \mathbb{C}^+ : |s| \geq \rho\}$. For $G_{cl} = (I + GK)^{-1}$ this follows from condition (iii), for terms containing $K$, $G$ we invoke (i). This proves $T \in \mathbf{H}^\infty$, hence the closed-loop is externally stable.

7) Since by our standing assumption controllers $K \in \mathscr{K}$ are admissible for $G$, the closed-loop system is well-posed. Since both $G$ and $K$ are exponentially stabilizable and detectable by (iv), and since the cross product $G \times K$ preserves these properties, the closed-loop system $T(G,K)$ is exponentially stabilizable and detectable by Morris [29, Thm. 6.1]. Therefore, by Lemma 1, exponential stability of the closed loop follows from its externally stability, which we proved in 5). That completes the proof. $\qquad\square$

**Remark 1.** The authors of [34] propose $h(s) = \left(\frac{s^2+s}{s^2+s+1}\right)^p$ for a pole of $F$ of order $p$ at 0, and similar expressions apply to poles off the origin. Multiplying with $h$ assures that the modified Nyquist curve (6) does not escape to infinity, as would be the case for more standard Nyquist curves with small $\varepsilon$-half circle indentations around open loop poles on $j\mathbb{R}$. This is favorable for its approximation by a polygon. The case occurs for instance in PID-control, see sections 7 - 10.

**Remark 2.** As simple an example as $G(s) = (s-1)^{-1}$ and $K(s) = (s-1)/(s+1)$ gives $n_p = 1$ and $\tilde{n}_p = 0$, which shows that pole zero cancellations may indeed occur. Our argument shows that in the case of a pole zero cancellation condition (6) is simply never satisfied. So our test cannot go wrong in that case.

**Remark 3.** In many applications the spectrum of the infinitesimal generator $A$ may be separated into two parts $\sigma_\pm(A)$ by a closed curve $\Gamma$ with $\sigma_-(A)$ lying outside $\Gamma$, and $\sigma_+(A)$ lying inside $\Gamma$ such that $\sigma_-(A) \subset \mathbb{C}^-$ and $\sigma_+(A)$ is discrete, hence finite. Then $\Sigma$ may be represented as the cross product of two systems $\Sigma_- \times \Sigma_+$, where $\Sigma_-$ is exponentially stable and $\Sigma_+$ is finite-dimensional. In that case hypothesis (iv) has only to be checked for $G_+$ (and $K$), which reduces to standard finite-dimensional tests like the Hautus test [20].

**Remark 4.** The interest in proving exponential stability of the closed loop lies of course in the well-known fact that it is preserved under linearization: If the Fréchet linearization about steady state of a nonlinear regulator $\mathbf{K}$ stabilizes the Fréchet linearization about steady-state of nonlinear system $\mathbf{G}$ exponentially, then $\mathbf{K}$ stabilizes $\mathbf{G}$ locally exponentially around that steady state. For infinite dimensional systems this is a consequence of Zwart [35].

**Remark 5.** As we shall see in the sequel, Theorem 1 gives key information for our algorithmic approach. As soon as unstable poles of $f$ and $F$ are the same for the initial stabilizing controller $K_0$, so that the winding number has the correct value, our method will only have to assure that the winding number does not change as the controller $K$ is updated during optimization. This guarantees that no unstable cancellations occur during optimization. Our technique to avoid changes of the winding number uses a barrier function and will be discussed in section 4.

## 3. Sampled Nyquist test with certificate

In this section we examine how the Nyquist test (6) is implemented. Writing $f$ instead of $\tilde f$, we seek $N \in \mathbb{N}$ and frequencies $\omega_0 = 0 < \omega_1 < \cdots < \omega_N = \infty$ such that the closed polygon $P_f = \{f(-j\omega_N), \dots, f(0), f(j\omega_1), \dots, f(j\omega_N)\}$ has the same winding number as the Nyquist curve $\{f(j\omega) : \omega \in \mathbb{R} \cup \{\infty\}\}$. Let $P_f(j\omega)$ denote the linearly interpolated function associated with the polygon, and for any $g$ let $\Delta_{[\omega', \omega'']} \arg g$ denote the change of argument of $g$ along the section $[j\omega', j\omega'']$ of $j\mathbb{R}$. Suppose $f(j\omega) \neq 0$ and $P_f(j\omega) \neq 0$ for all $\omega \in \mathbb{R} \cup \{\infty\}$, then, with the convention $\omega_{-i} = -\omega_i$ we have

$$\mathrm{ind}(f(j\mathbb{R}), 0) = -\frac{1}{2\pi} \sum_{i=-N}^{N-1} \Delta_{[\omega_i, \omega_{i+1}]} \arg f$$

and similarly

$$\mathrm{ind}(P_f, 0) = -\frac{1}{2\pi} \sum_{i=-N}^{N-1} \Delta_{[\omega_i, \omega_{i+1}]} \arg P_f$$
$$= -\frac{1}{2\pi} \sum_{i=-N}^{N-1} \arg[f(j\omega_{i+1})/f(j\omega_i)],$$

the last expression being computable. We now assure that these two winding numbers agree, which is true if the nodes $\omega_i$ are chosen such that, for every $i$,

$$\Delta_{[\omega_i, \omega_{i+1}]} \arg f = \arg[f(j\omega_{i+1})/f(j\omega_i)]. \tag{7}$$

Geometrically (7) means the closed curve $\gamma_i$ obtained by concatenating the segment $[f(j\omega_{i+1}), f(j\omega_i)]$ with the piece $f([j\omega_i, j\omega_{i+1}])$ of the Nyquist contour, does not encircle the origin. If $f(s)$ is available analytically we may, after fixing a small threshold $\delta > 0$, construct the $\omega_i$ through the recursion:

$$\omega_{i+1} = \sup\left\{\omega : \delta + \mathrm{Re}\int_{\omega_i}^{\omega}\frac{f'(j\omega)}{f(j\omega)}d\omega \le \arg\left[\frac{f(j\omega)}{f(j\omega_i)}\right]\right\}. \tag{8}$$

Alternatively, since $\arg[f(j\omega)/f(j\omega_i)] < \pi$, we may use the following slightly more conservative construction

$$\omega_{i+1} = \sup\left\{\omega : \delta + \mathrm{Re}\int_{\omega_i}^{\omega}\frac{f'(j\omega)}{f(j\omega)}d\omega \le \pi\right\}. \tag{9}$$

A third possibility to ensure (7) uses a bound on $f'$. Call $L[\cdot, \cdot]$ a first-order bound of $f$ if $L[\omega^-, \omega^+] \ge |f'(j\omega)|$ for every $\omega \in [\omega^-, \omega^+]$. Then we have the following simple test.

**Lemma 2.** *Let $\omega_i, \omega_{i+1}$ denote two consecutive nodes in the polygon $P_f$ not passing through $0$, and suppose*

$$L[\omega_i, \omega_{i+1}](\omega_{i+1} - \omega_i) < |f(j\omega_i)| + |f(j\omega_{i+1})|. \tag{10}$$

*Then condition (7) is satisfied.*

*Proof.* Assume on the contrary that the curve $\gamma_i$ in (7) encircles the origin. Let $\ell$ be the length of the curved part $\tilde\gamma_i = f([j\omega_i, j\omega_{i+1}])$ of $\gamma_i$. The projection line of 0 onto the segment $[f(j\omega_{i+1}), f(j\omega_i)]$ meets the segment at $P_f(j\omega^*)$, $\omega^* \in [\omega_i, \omega_{i+1}]$. But $\tilde\gamma_i$ has to cross this line at some point $p \notin [0, P_f(j\omega^*)]$ with $0 \in [p, P_f(j\omega^*)]$, so going from $f(j\omega_i)$ to $p$, $\tilde\gamma_i$ has length $\ge |f(j\omega_i)|$. Similarly, between $p$ and $f(j\omega_{i+1})$ the length of $\tilde\gamma_i$ is at least $|f(j\omega_{i+1})|$. Altogether, the length $\ell$ of $\tilde\gamma_i$ exceeds $|f(j\omega_i)| + |f(j\omega_{i+1})|$. But $\ell = \int_{j\omega_i}^{j\omega_{i+1}} |f'(z)|dz = \int_0^1 |f'(j\omega_i + t(j\omega_{i+1} - j\omega_i))|(\omega_{i+1} - \omega_i)dt \le L[\omega_i, \omega_{i+1}](\omega_{i+1} - \omega_i) < |f(j\omega_i)| + |f(j\omega_{i+1})|$ by hypothesis (10), a contradiction. $\square$

**Remark 6.** Consider the case where $G, K$ are stable so that $f$ is holomorphic on a domain containing $\overline{\mathbb{C}}^+$. Assume that $f$ is even holomorphic on $\mathbb{C}_{-\alpha}^+$ for some $\alpha > 0$, which is often the case, e.g. when $G$ is sectorial [21]. By (i) find $\beta > 0$ such that $f(\mathbb{C}_{-\alpha}^+) \subset \mathbb{C}_{-\beta}^+$, and put $\tilde f(s) = f(s - \alpha) + \beta$, then $\tilde f : \mathbb{C}^+ \to \mathbb{C}^+$. If $\Gamma = f(\gamma)$, $\gamma = j\mathbb{R}$, is the Nyquist curve, then $\tilde\Gamma = \tilde f(\tilde\gamma) = \Gamma + \beta$, where $\tilde\gamma = \alpha + j\mathbb{R}$, and in the place of $\mathrm{ind}(\Gamma, 0)$ we are now interested in $\mathrm{ind}(\tilde\Gamma, \beta)$. Put $\sigma(z) = \frac{1+z}{1-z}$ and $\tau(s) = \frac{s-1}{s+1}$, then $\tau = \sigma^{-1}$, and $\phi := \tau \circ \tilde f \circ \sigma$ maps the unit disk $\mathbb{D}$ to itself. Now $\gamma_0 = \tau(\tilde\gamma) = \{\frac{\alpha - 1 + j\omega}{\alpha + 1 + j\omega} : \omega \in \mathbb{R}\}$ is the circle with center $\frac{\alpha}{\alpha+1}$ and radius $\frac{1}{\alpha+1}$, the analogue of the Nyquist curve is $\Gamma_0 = \phi(\gamma_0) \subset \mathbb{D}$, and we are interested in $\mathrm{ind}(\Gamma_0, \tau(\beta))$. By the Schwarz-Pick theorem we have $|\phi'(z)| \le \frac{1 - |\phi(z)|^2}{1 - |z|^2}$ for $z \in \mathbb{D}$, hence $|\tilde f'(s)| \le \frac{2}{|1 + \phi(\tau(s))|^2}\frac{1 - |\phi(\tau(s))|^2}{1 - |\tau(s)|^2}\frac{2}{|1 + s|^2}$. We have to evaluate $f'(j\omega) = \tilde f'(\alpha + j\omega)$. Since $\tilde f(\alpha + j\omega) = f(j\omega)$ has a limit $\neq 0$ as $\omega \to \infty$, we have $\phi(\tau(\alpha + j\omega)) \not\to -1$, so that the term

$$\frac{1 - |\phi(\tau(\alpha + j\omega))|^2}{|1 + \phi(\tau(\alpha + j\omega))|^2}$$

remains bounded. But $\frac{1}{1 - |\tau(\alpha + j\omega)|^2}\frac{1}{|1 + \alpha + j\omega|^2} = \frac{1}{4\alpha}$, hence $\tilde f'$ is bounded on $\tilde\gamma = \alpha + j\mathbb{R}$, and so $f'$ is bounded on $j\mathbb{R}$. Using this, it also follows that $\phi'$ is bounded on $\gamma_0$. Going back with this information, we find that $|f'(j\omega)| = |\tilde f'(\alpha + j\omega)| \le C\omega^{-2}$ for some computable $C > 0$. This shows that for large $\omega$ the next frequency $\omega^+$ for the Nyquist sampling in the test (10) is of the order $\omega \sim \omega + C^{-1}\omega^2|f(j\infty)|$, which explains the extremely fast convergence in algorithm 1. It also follows that the first-order bound $L[\cdot, \cdot]$ is of the form $L[\omega, \omega^+] = C\omega^{-2}$. $\blacksquare$

Once $P_f$ is constructed, $\mathrm{ind}(P_f, 0)$ is computed by the ray-crossing algorithm: Fix a ray at the origin not passing through any of the nodes of $P_f$, and count in a straightforward way the number of signed crossings of that ray by the polygon $P_f$. The overall Nyquist procedure is presented in algorithm 1.

The following observation gives a justification of our approach.

**Algorithm 1.** Grid construction and Nyquist stability test

---

**Parameters:** $\delta > 0$, $\Theta > 1$.

▷ **Step 1** (`Initialize`). Choose $\omega_0 = 0$ and $\omega_1 > 0$ such that (9), respectively (10), is satisfied.

▷ **Step 2** (`Extrapolate`). Having constructed $\omega_0 < \cdots < \omega_i$, put $\omega^\sharp = \Theta(\omega_i - \omega_{i-1}) + \omega_i$. If (9), respectively (10), is satisfied on $[\omega_i, \omega^\sharp]$, then put $\omega_{i+1} = \omega^\sharp$, otherwise use backtracking to find $\omega_{i+1} \in (\omega_i, \omega^\sharp)$ such that (9), respectively (10), holds.

▷ **Step 3.** If $\omega_{i+1} < \infty$ loop on with step 2, otherwise obtain Nyquist grid $\Omega_{\text{nyq}}$ and goto step 4.

▷ **Step 4** (`Compute winding number`). Choose ray starting from 0 which avoids all $f(j\omega_i)$. Then count signed ray crossings of $P_f$ to obtain $\text{ind}(P_f, 0)$.

---

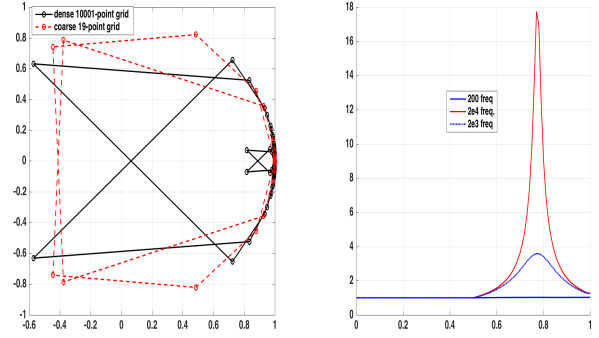

**Figure 1.** Comparison of logspace and adapted grid $\Omega_{\text{nyq}}$ for Nyquist (left). Cross section of $\|S\|_\infty$ on segment $[K, K^+]$ (right) for different grids.

**Theorem 2.** *Suppose for a given $K \in \mathscr{K}$ the integrals in (9) are computed formally to construct $P_f$, or a first-order bound $L[\cdot, \cdot]$ for $f$ is used to construct $P_f$ according to rule (10). Then the computation of the winding number is exact, i.e., satisfies $\text{ind}(f, 0) = \text{ind}(P_f, 0)$. In particular, if $\text{ind}(P_f, 0) = -n_p$, then $K$ is* certified *closed-loop stabilizing.* ∎

**Example 1.** Consider study 'DLR1' from the CompLeib collection [19], an open-loop stable rational system $G(s)$ with 10 modes, 2 control inputs and 2 measurements. All modes are badly damped and manifest themselves as sharp peaks in the frequency response with damping not exceeding 5e-3. With $K = [1, -1; -1, 1]$ the system is stable in closed loop, but when moving to $K^+ = [-1\ 1; 1\ -1]$, the closed-loop has two unstable modes $0.0041 \pm j0.9951$. Since the number of open-loop poles is $n_p = 0$, we expect the winding number 2 for $f = \det(I + GK^+)$ in (6).

We compute the index via ray-crossing of $P_f$ first on a dense grid $[0, \text{logspace}(-3, 3, 1000)]$, where we get the incorrect value 0. In contrast the grid of algorithm 1 needs only 19 frequencies with (9), and 27 with (10), yet delivers the correct winding number 2, which differs from $n_p = 0$, indicating the arrival of two unstable modes in closed loop. Fig. 1 (left) shows $P_f$ for $f = \det(I + GK^+)$ on the two grids. ∎

**Remark 7.** In Example 1 the closed-loop sensitivity $\|S(K + t(K^+ - K))\|_\infty$ has a bump at $t^* = .78$ on the segment $[K, K^+]$; see Fig. 1 (right). This is where instability occurs. The bump is more or less articulated depending on the frequency grid. This means that $\|S\|_\infty$ serves as a barrier against instability, but not always a reliable one, due to the fact that values re-descend as $t$ crosses $t^*$ and approaches 1.

## 4. Optimization method

In this section we present our algorithm for program (2). Let controllers $K \in \mathscr{K}$ be parametrized as $K(\mathbf{x})$ for some vector $\mathbf{x} \in \mathbb{R}^n$ of tunable parameters, and suppose the transfer functions $G(s)$ and $P(s)$ of system and plant are discretized on a sufficiently fine grid $\Omega_{\text{opt}} = \{\omega_0, \ldots, \omega_N\}$ for optimization. Then the closed-loop $H_\infty$-performance to be minimized is $h(\mathbf{x}) = \max_{\nu=0,\ldots,N} \overline{\sigma}(T_{zw}(P_w(j\omega_\nu), K(j\omega_\nu, \mathbf{x})))$. As square root of a maximum eigenvalue function, $h$ is locally Lipschitz, but nonsmooth and nonconvex.

Since we do not wish the Nyquist curve $f = \det(I + GK)$ to change its winding number as we update our controller $K(\mathbf{x})$ during optimization, we have to hinder $f$ from crossing 0. Using the sensitivity function $S = (I + GK)^{-1}$, this can be pursued by a constraint $\|S\|_\infty^{-1} \geq r^{-1}$ on the modulus margin, where $r > 0$ is some threshold. In discretized form this is a constraint

$$s(\mathbf{x}) := \max_{\nu=0,\ldots,N} \overline{\sigma}\left[(I + G_w(j\omega_\nu)K(j\omega_\nu, \mathbf{x}))^{-1}\right] \leq r. \quad (11)$$

Adding constraint (11) to program (2) gives the new cast

$$
\begin{array}{ll}
\text{minimize} & \max_{\nu=1,\ldots,N} \overline{\sigma}(T_{zw}(P(j\omega_\nu), K(j\omega_\nu, x))) \\
\text{subject to} & K \text{ stabilizes } G \text{ exponentially} \\
& (11) \text{ and } K \in \mathscr{K}
\end{array}
\quad (12)
$$

where the parameter $r$ is used to prevent crossings of the Nyquist curve, but serves also to improve robustness of the closed-loop system. Different robustness constraints against dynamic uncertainties could as well be included in the design. We refer the reader to [36, chap. 8] for a discussion.

In algorithm 2, instead of calibrating $r$, we use a dual approach, where we minimize the unconstrained function $g(\mathbf{x}) = \max\{h(\mathbf{x}), as(\mathbf{x})\}$ for some small penalty $a > 0$ over the parameter space $\mathbf{x} \in \mathbb{R}^n$. Since a maximum of $H_\infty$-norms is again an $H_\infty$-norm, Clarke subgradients of $g$ may be computed by the method of [1]. We now apply algorithm 2 to minimize $g$, that is, to solve program (2).

During the following we comment on the salient features of this scheme.

**Remark 8.** The *primary descent step* $\mathbf{x}^+$ of step 3 is simply the standard step of the nonsmooth trust-region method [4].

---

**Algorithm 2.** Non-smooth optimization for (2)

▷ **Step 1** (`Initialize`). Find initial stabilizing controller $K(\mathbf{x}^0)$, put counter $j = 0$, and determine number $n_p$ of open-loop rhp poles for Nyquist test. Initialize trust-region radius as $R_1 > 0$.

▷ **Step 2** (`Local model`). Given current iterate $\mathbf{x}^j$ at counter $j$, compute a local polyhedral model $\phi(\cdot, \mathbf{x}^j)$ of $g$ at $\mathbf{x}^j$.

▷ **Step 3** (`Primary descent`). Starting with trust-region radius $R_j$ and model $\phi$, use trust-region update mechanism in tandem with local model update to generate a primary descent step $\mathbf{x}^+$ for $g$. Procedure ends with new trust-region radius $R^+$, and new local model $\phi^+(\cdot, \mathbf{x}^j)$.

▷ **Step 4** (`Nyquist test`). Use Nyquist test in algorithm 1 to check whether $K(\mathbf{x}^+)$ is closed-loop stabilizing. If this is the case (i.e., $\mathrm{ind}(P_f, 0) = -n_p$), then put $\mathbf{x}^{j+1} = \mathbf{x}^+$ and $R_{j+1} = R^+$, increase counter $j$, and loop on with step 2. In case of instability ($\mathrm{ind}(P_f, 0) \neq -n_p$), goto step 5.

▷ **Step 5** (`Stability safeguard`). Reject descent step $\mathbf{x}^+$, reduce trust-region radius to $R^{++} = R^+/2$, and add a repelling cutting plane to the local model $\phi^+$ to obtain $\phi^{++}$. Then go back to step 3 with initial information $R^{++}$, $\phi^{++}$ instead of $R_j$ and $\phi$.

---

Here $\mathbf{x}^+$ gives sufficient decrease of $g$, and would normally be accepted as the next serious iterate. The trouble is that $\mathbf{x}^+$ may lead to a destabilizing controller $K(\mathbf{x}^+)$.

The difficulty is explained as follows. In the majority of cases the function $s : \mathbf{x} \mapsto \|S(\mathbf{x})\|_\infty$ has a barrier effect as it makes $K(\mathbf{x})$ approach the boundary of stability from inside (see Fig. 1 right). But in contrast with classical barrier functions like the log-barrier in interior point methods, $s(\mathbf{x})$ takes on finite values behind the barrier and outside the domain of stability. This means it cannot be fully relied on to enforce stability, as seen in Example 1. This is why it is used in tandem with the Nyquist test.

**Remark 9.** In the original approach [1, 37] to nonsmooth $H_\infty$-synthesis the closed-loop system matrix $A(K)$ is available, so that closed-loop stability can be implemented using the spectral abscissa: a constraint $\alpha(A(K(\mathbf{x}))) \leq -\varepsilon$ is added, which not only serves to *recognize* instability, but also allows to *repel* steps from becoming unstable. In contrast, our Nyquist test allows to detect instability, but since the winding number is a discrete quantity, it cannot be used as a constraint to generate the repelling effect. The latter is implemented through the barrier property of the sensitivity function (11). Backtracking from the unstable $\mathbf{x}^+$ toward the stable $\mathbf{x}^j$, we locate an intermediate stable value $\mathbf{x}_t = t\mathbf{x}^j + (1-t)\mathbf{x}^+$, for which $s(\mathbf{x}_t)$ is large but $K(\mathbf{x}_t)$ is still stabilizing. Then we generate a cutting plane of $s(\cdot)$ at $\mathbf{x}_t$, which we add to the model

$\phi^+$. Ideally, this plane is relatively steep and therefore builds the desired barrier effect into the improved model $\phi^{++}$.

**Remark 10.** Step 1 requires that $G$ can be stabilized by a structured controller $K_0 \in \mathscr{K}$. This is a stronger hypothesis than in (iv). Even for finite-dimensional systems it is generally difficult to *decide* whether a stabilizing controller of a given structure $\mathscr{K}$ exists. The problem is known to be NP-hard for static, reduced fixed-order, or PID controllers. However, this is a worst-case result which is somewhat in contrast with the fact that practical systems are usually easy to stabilize even with structured laws.

Note that when $G$ is already stable, an initial stabilizing controller is obtained by the small gain condition $\|K(s)\|_\infty < 1/\|G(s)\|_\infty$. This guarantees $T \in \mathbf{H}^\infty$ in (5), and then internal stability under the hypotheses of Theorem 1.

**Remark 11.** Convergence analysis of the trust-region algorithm is outside the scope of this work and can be based on [4]. Note that due to nonsmoothness special care has to be taken, as the standard trust-region scheme based on the Cauchy point fails. The success of the method hinges on building a good polyhedral model of the objective at the current iterate based on cutting planes. We refer to [1], where this is discussed.

## 5. Sampling for synthesis with certificate

As we have seen, Nyquist stability (6) can be based on the relatively coarse grid $\Omega_{\mathrm{nyq}}$ of algorithm 1. This typically requires significantly less then 100 nodes for most plants, but $\Omega_{\mathrm{nyq}}$ must at each step be re-adapted to the candidate controller $K(\mathbf{x})$, because the tunable parameters $\mathbf{x}$ move during optimization. In contrast, the grid for optimization $\Omega_{\mathrm{opt}}$ used in algorithm 2 has to be of finer scale, but it remains invariant during optimization, as updating it would change problem (2).

An initial stabilizing controller $K(\mathbf{x}_0)$ can be used to build an initial grid $\Omega_{\mathrm{opt}}$, but we have to be aware that the controller $K(\mathbf{x})$, which varies in the course of optimization, may by itself develop resonant modes, which may render the original sampling $\omega_i \in \Omega_{\mathrm{opt}}$ inappropriate. In order to prevent this phenomenon, it is cautious to optimize over classes $\mathscr{K}$ of *stable* controllers, and to put *constraints on the damping* of the controller modes, confining them to a conical region in $\mathbb{C}^-$. For real-rational controllers with explicit state-space realization (3) such constraints are readily implemented and added in (2). Stability of $K$ translates into $\alpha(A_K(\mathbf{x})) \leq -\varepsilon$ for some threshold $\varepsilon > 0$, and similarly, the mode damping requirement becomes a constraint $\alpha_\triangleright(A_K(\mathbf{x})) \leq -r$ as soon as we define a conical analogue of the spectral abscissa via

$$\alpha_\triangleright(A) = \max\{\mathrm{Re}(\lambda)/|\lambda| : \lambda \text{ eigenvalue of } A\},$$

and where $r$ accounts for the aperture of the conical section. Subgradients of $\alpha$ and $\alpha_\triangleright$ are computed as in [38, 39].

Even with these precautions, upgrading $\Omega_{\text{opt}}$ may become necessary, and we now discuss a method to adapt a grid $\Omega_{\text{opt}}$ to a candidate controller $K \in \mathscr{K}$.

**Lemma 3.** *Let $\phi : \mathbb{R} \to \mathbb{R}_+$ be of class $C^1$, and let $L[\cdot,\cdot]$ be a first-order bound for $\phi$. Let $\omega_i, \omega_{i+1}$ be two consecutive nodes of a piecewise linear interpolation $P_\phi$ of $\phi$ such that $\gamma^* \geq \max\{\phi(\omega_i), \phi(\omega_{i+1})\}$ and for a given tolerance $\vartheta > 0$,*

$$L[\omega_i, \omega_{i+1}](\omega_{i+1} - \omega_i) < 2\gamma^* + 2\vartheta - \phi(\omega_i) - \phi(\omega_{i+1}). \quad (13)$$

*Then $\phi(\omega) < \gamma^* + \vartheta$ for every $\omega \in [\omega_i, \omega_{i+1}]$.*

*Proof.* Suppose on the contrary that there exists $\omega^* \in [\omega_i, \omega_{i+1}]$ such that $\phi(\omega^*) \geq \gamma^* + \vartheta$. Then the polygon connecting $\phi(\omega_i), \phi(\omega^*), \phi(\omega_{i+1})$ has length $\geq L$, $L = \sqrt{A^2 + (\omega^* - \omega_i)^2} + \sqrt{B^2 + (\omega_{i+1} - \omega^*)^2}$, where $A = \gamma^* + \vartheta - \phi(\omega_i)$ and $B = \gamma^* + \vartheta - \phi(\omega_{i+1})$. We have $L \geq \ell = \sqrt{(A+B)^2 + (\omega_{i+1} - \omega_i)^2}$, the minimum being attained at $\omega^* = \frac{\omega_i B + \omega_{i+1} A}{A+B}$. But the curve $\{(\omega, \phi(\omega)) : \omega \in [\omega_i, \omega_{i+1}]\}$ has length

$$\mathscr{L} = \int_{\omega_i}^{\omega_{i+1}} \sqrt{1 + \phi'(\omega)^2}\, d\omega \leq \sqrt{1 + L[\omega_i, \omega_{i+1}]^2}(\omega_{i+1} - \omega_i),$$

and $\mathscr{L} \geq L \geq \ell$, so we get the estimate $\sqrt{1 + L[\omega_i, \omega_{i+1}]^2} \geq \sqrt{(A+B)^2/(\omega_{i+1} - \omega_i)^2 + 1}$, which contradicts (13). $\qquad\square$

---

**Algorithm 3.** Infinite-dimensional $H_\infty$-synthesis

**Parameters:** Tolerance $\vartheta > 0$.

▷ **Step 1** (`Grid for optimization`). Use initially stabilizing controller $K_0 \in \mathscr{K}$ and first-order bound condition (13) for function $\phi(\omega) = \overline{\sigma}(T_{wz}(P(j\omega), K_0(j\omega)))$ to construct grid $\Omega_{\text{opt}}$.

▷ **Step 2** (`Optimize`). Using algorithm 2, compute optimal controller $K^* \in \mathscr{K}$ on grid $\Omega_{\text{opt}}$ with value $\gamma^*$.

▷ **Step 3** (`Refined grid`). Use first-order bound $L[\cdot,\cdot]$ for $\phi(\omega) = \overline{\sigma}(T_{wz}(P(j\omega), K^*(j\omega)))$ to check whether grid $\Omega_{\text{opt}}$ satisfies (13). If not add nodes to assure this and obtain verification grid $\Omega_{\text{ver}}$ with this property.

▷ **Step 4** (`Verify`). Check $\gamma^* \geq \max_{\Omega_{\text{ver}}} \overline{\sigma}(T_{wz}(P, K^*)) - \vartheta$. If this is the case quit successfully, otherwise replace $\Omega_{\text{opt}}$ by $\Omega_{\text{opt}} \cup \Omega_{\text{ver}}$ and go back to step 2.

---

As we shall see in the sequel, $\vartheta > 0$ serves as the tolerance within which we are able to know the value of the infinite-dimensional (un-sampled) $H_\infty$-norm $\|T_{wz}(P, K)\|_\infty$. In order to derive this, we have to apply the test (13) to the performance function $\phi(\omega) = \overline{\sigma}(T_{wz}(P(j\omega), K(j\omega)))$, and for that we have to analyze its differentiability. Consider the one-parameter family of symmetric matrices

$$\omega \mapsto \mathscr{M}(\omega) = T_{wz}(P(j\omega), K(j\omega))^H\, T_{wz}(P(j\omega), K(j\omega)),$$

then by [40, Theorem 6.1] the eigenvalues $\lambda_\nu(\omega)$ of $\mathscr{M}(\omega)$ are real analytic functions, hence $\phi^2$ is a finite maximum of

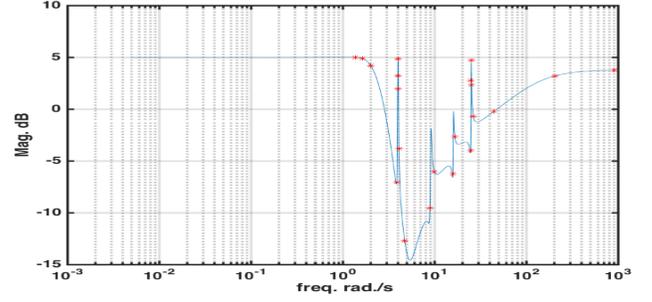

**Figure 2.** Verification grid $\Omega_{\text{ver}}$ via (13), with $\gamma^* = 1.78$ and $\vartheta = 10^{-2}$. As expected flat parts need few grid points $\omega_i$, whereas resonances are perfectly captured.

real analytic eigenvalue functions, and since $\phi > 0$, we deduce that $\phi$ as well is a finite maximum of real-analytic functions. What is even more important is the following:

**Lemma 4.** [3, Theorem 2.3] $\phi$ *has only finitely many points of non-smoothness, and is of class $C^2$ at peak frequencies.* $\quad\square$

In consequence, there exists $\vartheta_0 > 0$ such that $\phi$ is of class $C^2$ on $\{\omega : \phi(\omega) > \|T_{wz}(P, K)\|_\infty - \vartheta_0\}$. This leads to

**Theorem 3.** *If $0 < \vartheta \leq \vartheta_0$ and a first-order bound $L[\cdot,\cdot]$ for $\phi = \overline{\sigma}(T_{wz}(P, K^*))$ in tandem with rule (13) is used in step 4 of algorithm 3, then the gain $\gamma^*$ achieved by $K^*$ is certified to satisfy*

$$\gamma^* \geq \|T_{wz}(P, K^*)\|_\infty - \vartheta. \quad (14)$$

$\blacksquare$

**Remark 12.** In practice it is usually sufficient to generate a numerical upper bound $L[\cdot,\cdot]$ using a finite-difference approximation $\phi'(\omega) \approx (\phi(\omega^+) - \phi(\omega^-))/(\omega^+ - \omega^-)$. In our testing this gives excellent results and leads to moderately sized grids $\Omega_{\text{opt}}$ and $\Omega_{\text{ver}}$. Fig. 2 gives a typical case.

**Remark 13.** To generate the optimization grid $\Omega_{\text{opt}}$ we apply (13) with $\gamma^* = \max\{\phi(\omega_i), \phi(\omega_{i+1})\}$ on each interval $[\omega_i, \omega_{i+1}]$. When it comes to just certifying the optimal value $h(\mathbf{x}^*) = \|T_{wz}(K(\mathbf{x}^*))\|_{\infty,d} = \phi(\omega^*)$ in step 4 of algorithm 3, we can construct an even coarser grid by applying (13) with $\gamma^* = \phi(\omega^*)$ on every $[\omega_i, \omega_{i+1}]$. Here our grid turns out sparse at frequencies $\phi(\omega) \ll \phi(\omega^*)$, while resonances are still accurately captured (see Fig. 2 for an illustration). We call this a verification grid $\Omega_{\text{ver}}$. The outlined method to construct $\Omega_{\text{opt}}$, and to complete it in step 4 by adding elements of a verification grid $\Omega_{\text{ver}}$, is well-suited to discretize the controller design problem (1). Discretization at that level avoids the pitfalls in system reduction and identification techniques.

We can further exploit Lemma 3 to obtain information on how close the values $\gamma^*$ of (2), and $\gamma_\infty$ of (1), are. Writing

$h(\mathbf{x}) = \|T_{wz}(K(\mathbf{x}))\|_{\infty,d}$ for the discrete $H_\infty$-norm of (2) on $\Omega_{\text{opt}}$, and $h_\infty(\mathbf{x}) = \|T_{wz}(K(\mathbf{x}))\|_\infty$ for the true $H_\infty$-norm in (1), we have the following:

**Corollary 1.** *Let $\mathbf{x}_\infty$ be a local minimum of the infinite-dimensional $H_\infty$-program with value $\gamma_\infty$, and $\mathbf{x}^*$ a local minimum of (2) with value $\gamma^*$. Suppose a first-order bound in tandem with rule* (13) *is used in step 3 of algorithm 3. Then if $\mathbf{x}^*$, $\mathbf{x}_\infty$ are within neighborhoods of local optimality of each other, we have $h(\mathbf{x}_\infty) \geq h(\mathbf{x}^*) \geq h_\infty(\mathbf{x}^*) - \vartheta \geq h_\infty(\mathbf{x}_\infty) - \vartheta \geq h(\mathbf{x}_\infty) - \vartheta$.*

*Proof.* Indeed, $h(\mathbf{x}_\infty) \geq h(\mathbf{x}^*)$ because $\mathbf{x}^*$ is a minimum of $h$ on a neighborhood $U(\mathbf{x}^*)$, and $\mathbf{x}_\infty \in U(\mathbf{x}^*)$ by hypothesis. Next $h(\mathbf{x}^*) \geq h_\infty(\mathbf{x}^*) - \vartheta$ by Lemma 3, because construction of the grid uses the bound $L[\cdot,\cdot]$ and rule (13). Next $h_\infty(\mathbf{x}^*) \geq h_\infty(\mathbf{x}_\infty)$, because $\mathbf{x}_\infty$ is a minimum of $h_\infty$ on a neighborhood $U(\mathbf{x}_\infty)$, and $\mathbf{x}^* \in U(\mathbf{x}_\infty)$ by hypothesis. The last inequality is satisfied because $h \leq h_\infty$. $\qquad\square$

This means comparable locally optimal values of the infinite dimensional $H_\infty$-program (1) and its approximation (2) differ by at most $\vartheta$, our apriori chosen tolerance. Since most of the time our algorithm finds even the global minimum of (2), this is a very useful information in practice, as the value $\gamma_\infty$ of a global solution of the infinite dimensional $H_\infty$-program is then known within the prior tolerance $\vartheta$.

The result of Theorem 3 could also be explained as follows. Suppose $\mathbf{x}^*$ is a local minimum of (2), i.e., $h(\mathbf{x}) \geq h(\mathbf{x}^*)$ for every $\mathbf{x}$ in some neighborhood $U$ of $\mathbf{x}^*$. We know that $h \leq h_\infty$, so the value $\gamma^* = h(\mathbf{x}^*)$ is *a priori* optimistic. Could it be overly optimistic (i.e. could it be way too low) and therefore misleading? The answer is no.

**Corollary 2.** *Let $\gamma_\infty$ be the best value of program* (1) *on $U$, that is $\gamma_\infty = \inf\{h_\infty(\mathbf{x}) : \mathbf{x} \in U \text{ admissible in } (1)\}$. Then $\gamma^* \leq \gamma_\infty \leq \gamma^* + \vartheta$.*

*Proof.* Since $h \leq h_\infty$ we have $\gamma^* = \inf_U h \leq \inf_U h_\infty = \gamma_\infty$. Fix $\varepsilon > 0$, then there exists $\mathbf{x}_\infty \in U$ such that $\gamma_\infty \geq h_\infty(\mathbf{x}_\infty) - \varepsilon$. By Theorem 3 we have $\gamma^* \geq h_\infty(\mathbf{x}^*) - \vartheta \geq \gamma_\infty - \vartheta \geq h_\infty(\mathbf{x}_\infty) - \vartheta - \varepsilon \geq h(\mathbf{x}_\infty) - \vartheta - \varepsilon \geq h(\mathbf{x}^*) - \vartheta - \varepsilon = \gamma^* - \vartheta - \varepsilon$, and since $\varepsilon$ is arbitrary, this implies $\gamma^* \geq \gamma_\infty - \vartheta \geq \gamma^* - \vartheta$. $\qquad\square$

## 6. Boundary and distributed PDE control

Developing Nyquist stability and $H_\infty$-optimization for well-posed transfer functions $G(s)$ has the advantage that a wide set of potential applications is covered. In this section we illustrate our strategy for distributed and boundary control of partial differential equations. Numerical tests are included in section 10. Following [25, 21, 41, 42], a boundary control problem may be represented in the abstract form

$$\Gamma : \begin{cases} \dot{x} &= Ax \\ Px &= u \\ y &= Cx \end{cases} \qquad (15)$$

with operators $A \in L(X, H)$, $P \in L(X, U)$, $C \in L(X, Y)$ on Hilbert spaces $X, H, U, Y$, where $X$ is dense in $H$ and $D(A) \subset D(P)$. The idea developed by Salomon [25] is now to represent $\Gamma$ by a well-posed system $\Sigma_\Gamma$ with input $u$ and output $y$, thereby making it amenable to techniques developed for this class. As in [25, 41] one lets $X_0 = X \cap \ker(P)$ and restricts $A$ to $X_0$ to generate the semi-group, while $C$ restricted to $X_0$ induces the output operator. Construction of a suitable control operator $B$ is more involved, and we refer to [25] and [21] for details.

The transfer function $G(s)$ of $\Gamma$ can be obtained by applying the Laplace transform, [41], which leads to a family $\Gamma_s$ of abstract elliptic boundary control problem

$$\Gamma_s : \begin{cases} sx(s) &= Ax(s) \\ Px(s) &= u(s) \\ y(s) &= Cx(s) \end{cases} \qquad (16)$$

The question is then how well-posedness of $G(s)$ and conditions (i) - (iv) can be verified.

For parabolic and hyperbolic PDEs well-posedness was first examined in [25]. A systematic study is Cheng and Morris [41], where it is shown that under natural hypotheses $\Sigma_\Gamma$ is well-posed iff $G(s)$ is bounded on some half-plane $\text{Re}(s) > \sigma$. This is beneficial for our present approach in so fas as we do not have to construct $\Sigma_\Gamma$ explicitly, and can concentrate on carrying out synthesis in the frequency domain. Computation of $G(s)$ may be based either on a formal or a numerical evaluation of (16) at a given $s$. The remaining issue is then to check condition (iv) in a given situation.

**Example 2.** We consider boundary control of heat flow in a one-dimensional medium

$$\Gamma : \quad x_t(\xi, t) - x_{\xi\xi}(\xi, t) = 0, \qquad 0 \leq \xi \leq 1, t \geq 0$$

with initial conditions $x(\xi, 0) = 0$ and Neumann boundary control

$$x_\xi(0, t) = 0, \; x_\xi(1, t) = u(t),$$

where $u(t)$ is the rate of heat flow into the medium at the end $\xi = 1$. As measurement we take $y(t) = x(\xi_0, t)$ at some position $0 \leq \xi_0 \leq 1$. Following [21, Example 4.3.12], the transfer function $G(s) = y(s)/u(s)$ is

$$G(s) = \frac{1}{s} + 2\sum_{\nu=1}^\infty \frac{(-1)^\nu \cos(\nu\pi\xi_0)}{s + \nu^2\pi^2}$$

from which we see that $G$ is strictly proper and meromorphic, but not stable due to the pole at 0. Note that a closed form for $G(s)$ is given in [43], and the results in [41] show that $G(s)$ is well-posed. For well-posedness of the Dirichlet case see e.g. [42].

Before we apply the Nyquist test, we have to check hypothesis (iv). To see that the system is stabilizable we take the state feedback law $u(t) = -\alpha x(1, t)$ with $\alpha = \sqrt{k}\tan\sqrt{k} >$ 0 for some $k \in (0, \frac{\pi}{2})$, then the state evolves as $x(\xi, t) = x(0,0)e^{-kt} \cos \sqrt{k}\xi$, which decays exponentially in $t$ uniformly over $\xi \in [0,1]$.

For detectability, we have to find a law $F : h(t) \mapsto v(\xi)h(t)$ such that $x_t = x_{\xi\xi} + v(\xi)x(\xi_0, t)$ with boundary conditions $x_\xi(0, t) = 0 = x_\xi(1, t)$ is stable, and that works similarly. Experiments with this example are included in Section 10. For the setup of boundary control problems see also [25, 43]. ∎

**Remark 14.** In a well-posed system $\Sigma$ exponential stabilizability in (iv) can be verified by the following condition: For every initial state $x_0$ there exists $u \in L^2([0, \infty), U)$ such that the solution of $\dot{x} = Ax + Bu$, $x(0) = x_0$, is in $L^2([0, \infty), X)$. In [44] this is referred to as optimizability, and by [45, Thm. 1.1] is equivalent to exponential stabilizability. For exponential detectability one can use the formally weaker but equivalent estimatability [45]. For bounded $B, C$ equivalence is shown in [21]. The advantage is that these open-loop conditions can be checked in the original problem (15). In the context of $\Gamma$, given $x_0$, we have to make sure that we can find an open-loop boundary control $u \in L^2([0, \infty), U)$ such that the solution of $\dot{x} = Ax$, $x(0) = x_0$, $Px = u$ is in $L^2([0, \infty), X)$. This was used in example 2. Constructing $\Sigma_\Gamma$ explicitly may then again be avoided. For detectability, the situation is similar.

In boundary control of several spatial dimensions input and output spaces are usually infinite-dimensional, so that in order to comply with our standing hypothesis $U \simeq \mathbb{R}^p$, we may have to select a finite set of boundary basis functions $\phi_1, \ldots, \phi_p$ and restrict the boundary control operator $P$ to controls $u \in U$ of the form $u(t, \xi) = \sum_{i=1}^p u_i(t)\phi_i(\xi)$, $\xi \in \partial\Omega$. Similarly, $Y \simeq \mathbb{R}^m$ is usually achieved by taking a finite set of measurements $y_k(t) = \int_\Omega c_k(\xi, t)x(\xi, t) \, d\xi$ over the spatial domain $\Omega$. For problems with one spatial dimension point measurements are also possible [42]. These discretizations do not affect the question of well-posedness. Systematic ways to get finite-dimensional approximations of infinite dimensional controllers are discussed in Morris [46].

**Remark 15.** There is a rich literature on state-feedback stabilizability of boundary and distributed control problems for PDEs. For parabolic equations, where the semi-group is analytic [25], the spectrum decomposition condition is satisfied, so stabilizability can be checked using the Hautus test for the finite-dimensional subsystem, see [47, 20]. This has been exploited for a variety of parabolic equations. For the Navier-Stokes equation see e.g. [48], for a parallel heat flow exchanger see [49], for an unstable heat equation see [50]. Analytic semigroups preserve their favorable structure in closed loop $A + BKC$ with unbounded $B$ and bounded $KC$, as follows from the result in [51].

In finite-dimensional systems $G(s)$ a minimal realization is automatically stabilizable and detectable, so external stabilization by feedback will at least render a minimal realization of the closed loop exponentially stable. This may be considered a license to work directly in the frequency domain.

**Remark 16.** In contrast, even though minimal realizations for infinite-dimensional well-posed systems exist [30, Sect. 9], their value is limited, as they are not automatically stabilizable nor detectable. Logemann [52] gives the example $G(s) = (s+1)^{-1}(s(1 - e^{-s}) + 1)^{-1} \in \mathbf{H}^\infty(\mathbb{C}^+)$, which maps $L^2$ into $L^2$, yet its minimal realization is not exponentially stable.

As we cannot count on minimal realizations to assure a version of condition (iv), we propose the following result, which gives at least a partial remedy in infinite dimensions.

**Theorem 4.** Suppose $G(s)$ is in $L^2_\sigma(U, Y)$ for some $\sigma \geq 0$, and extends meromorphically into $\mathrm{Re}(s) > -\alpha$ for some $\alpha > 0$. Suppose $G(s) - G(\infty) = O(s^{-r})$ for some $r > \frac{1}{2}$ as $s \to \infty$ on $\mathrm{Re}(s) > -\alpha$. Suppose $K^* \in \mathscr{H}$ is computed by algorithm 3, hence satisfies (6). Then $G(s)$ admits a well-posed state-space realization with regard to which the closed loop with $K^*$ is exponentially stable.

*Proof.* By Mossaheb [33] the strictly proper $G(s) - G(\infty)$ has coprime factorizations over $\mathbf{H}^\infty(\mathbb{C}^+_\alpha)$ due to the sufficiently rapid decay $O(s^{-r})$, and hence so has $G(s)$. Since $G$ is the frequency representation of an operator $\mathscr{G} \in TIC_\sigma(U, Y)$, it follows from [30, Thm. 8.4.1(ii)] that $G(s)$ has a jointly exponentially stabilizable and detectable well-posed realization. But now all the hypotheses of Theorem 1 are satisfied, hence the Nyquist test (6) assures that semi-group of the closed-loop $T(G, K^*)$ is exponentially stable. □

**Remark 17.** Still in the same vein, when it is known that $G(s)$ can be realized by a well-posed system $\Sigma$ whose generator $A$ satisfies the spectrum decomposition condition, then on decomposing $\Sigma$ into its exponentially stable part $\Sigma_-$ and its finite-dimensional part $\Sigma_+$, we know that on taking a finite-dimensional realization $\widetilde{\Sigma}_+$ of the finite-dimensional part, and on patching $\Sigma_-$ and $\widetilde{\Sigma}_+$ together, we can always get a reduced well-posed system $\widetilde{\Sigma}$ representing $G(s)$, which satisfies hypothesis (iv). Unless we are specifically interested in analyzing stabilizability of the given representation, we may therefore in these cases avoid the explicit construction of $\Sigma_\Gamma$ and work directly with $\Gamma$, or entirely in the frequency domain.

**Remark 18.** Not unexpectedly, hyperbolic boundary and distributed control problems are more mulish with regard to applicability of our method. Well-posedness of such systems was first studied in [25], and [53] shows well-posedness for problems in one spatial dimension. For two and more spatial dimensions a case-by-case study is needed.

As far as conditions (i) - (iv) in Theorem 1 are concerned, the primary difficulty is that hyperbolic systems may have infinitely many unstable open-loop poles, in which case the Nyquist test is clearly not directly applicable. In that case our method may still be used in the optimization phase if an initial stabilizing controller is found by some other method. Even when there are only finitely many unstable poles, a second difficulty arises when stable poles accumulate along the imaginary axis. This may foil properness of $G$ in (i), but more typical is that exponential stabilizability of $G$ in (iv) fails. The well-known example of Renardy [54] shows that this may even happen when stable poles accumulate along a line $\text{Re}(s) = -\alpha$ with $\alpha > 0$.

For hyperbolic systems a version of Theorem 1 based on the notion of strong stability is helpful. Since the transfer function of hyperbolic systems is as a rule meromorphic on a domain containing $\overline{\mathbb{C}^+}$, the following result is interesting:

**Proposition 1.** *Let $K^* \in \mathcal{K}$ be computed by algorithm 3. Suppose $G, K^*$ are meromorphic on domain containing $\overline{\mathbb{C}^+}$, let* (i)-(iii) *be satisfied, and replace* (iv) *by the weaker condition* (iv') $G, K^*$ *are strongly stabilizable and strongly detectable. Then the closed loop $T(G, K^*)$ is strongly stable, and approximate optimality* (14) *is achieved.*

*Proof.* As in the proof of Theorem 1 we use $h$ to remove a finite number of open-loop poles on $j\mathbb{R}$. Since $G, K^*$ are meromorphic on a domain containing $\overline{\mathbb{C}^+}$, we can carry out the reasoning in Theorem 1 which gives $T \in \mathbf{H}^\infty$. Now according to [30, Lemma 8.2.7] the closed loop is also strongly stabilizable and strongly detectable, and since the closed-loop is input-output stable, we conclude using [30, Thm. 8.2.11] that the closed loop is even strongly stable. □

**Example 3.** In [43] the authors consider a damped wave equation where $G$ is open-loop externally stable, but stable poles accumulate along $j\mathbb{R}$ so that the system is strongly stable but not exponentially stabilizable. In this situation our method can still be used if we accept strong stability of the closed loop as satisfactory. For similar examples see e.g. [55].

## 7. Application in process control

We apply our frequency-sampled $H_\infty$-synthesis method to control a continuous cooling crystallizer, shown schematically in Fig. 3. The process uses fines dissolution and product removal, and is governed by a population balance and a molar balance equation; see [56, 57]. The population balance is of the form

$$\frac{\partial n(L,t)}{\partial t} = -G(c(t))\frac{\partial n(L,t)}{\partial L} - \frac{q}{V}h_{fp}(L)n(L,t) \quad (17)$$

$$n(L,0) = n_0(L), \qquad n(0,t) = \frac{B(c(t))}{G(c(t))} \quad (18)$$

where $n(L,t)$ is the crystal size distribution (CSD), $c(t)$ is the solute concentration, and the classification functions specifying fines dissolution and product removal are $h_f(L) = R(1 - h(L - L_f))$, $h_p(L) = 1 + zh(L - L_p)$, $h_{fp} = h_f + h_p$, where $h$ is the unit-step function. The crystal growth and birth coefficients obey phenomenological laws

$$G(c) = k_g(c - c_s)^g, \qquad B(c) = k_b(c - c_s)^b.$$

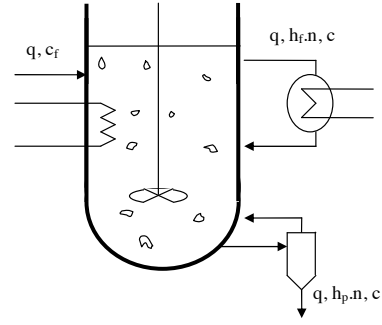

**Figure 3.** Continuous KCl-crystallizer with solute feed $c_f$, fines dissolution $h_f$, and product removal $h_p$. Solute concentration $c(t)$ is stabilized at steady-state by control of solute feed $c_f(t)$.

The molar balance is an integral-differential equation of the form

$$M\frac{dc}{dt} = \frac{q(\rho - Mc)}{V} + \frac{\rho - Mc}{\varepsilon}\frac{d\varepsilon}{dt} + \frac{qMc_f}{V\varepsilon} - \frac{q\rho}{V\varepsilon} - \frac{q\rho z\eta}{V\varepsilon}, \quad (19)$$

with initial condition $c(0) = c_0$, where

$$\varepsilon(t) = 1 - k_v \int_0^\infty n(L,t)L^3 dL, \qquad \eta(t) = k_v \int_{L_p}^\infty n(L,t)L^3 dL.$$

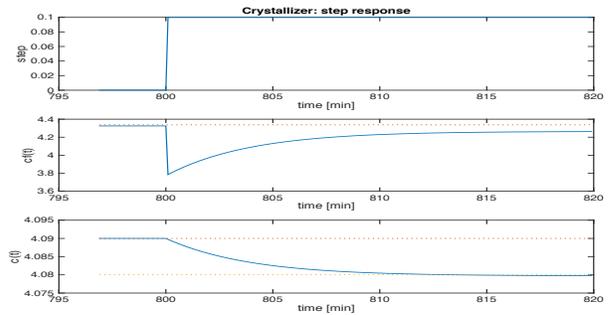

**Figure 4.** Open-loop step response of nonlinear model

The steady state equations lead to the explicit relationship

$$Mc_{fss} = \rho(1 + z\eta_{ss}) - (\rho - Mc_{ss})\varepsilon_{ss}$$

where

$$\varepsilon_{ss} = 1 - k_v \int_0^\infty n_{ss}(L)L^3 dL, \qquad \eta_{ss} = k_v \int_{L_p}^\infty n_{ss}(L)L^3 dL$$

with

$$n_{ss}(L) = \frac{B(c_{ss})}{G(c_{ss})}e^{-\frac{q}{VG(c_{ss})}H_{fp}(L)}, H_{fp}(L) = \int_0^L h_{fp}(\ell)d\ell,$$

and where our experiment uses $c_{ss} = 4.09$. Parameters are gathered in Table 1. The control input is solute feed concentration $c_f(t)$, the measured output is molar concentration $c(t)$.

Crystallizer data

| | | | |
|---|---|---|---|
| feed rate | $q$ | 0.05 | $\ell/min$ |
| total volume | $V$ | 10.5 | $\ell$ |
| fines removal size | $L_f$ | 0.2 | $mm$ |
| product removal size | $L_p$ | 1.0 | $mm$ |
| fines removal rate | $R$ | 5.0 | $--$ |
| product removal rate | $z$ | 2.0 | $--$ |
| growth rate constant | $k_g$ | 0.0305 | $mm\ell/min \cdot mol$ |
| growth rate exponent | $g$ | 1 | $-$ |
| nucleation rate | $k_b$ | 8.36e9 | $\ell^3/min \cdot mol^4$ |
| nucleation rate exponent | $b$ | 4 | $-$ |
| crystal density | $\rho$ | 1989 | $g/\ell$ |
| molar mass | $M$ | 74.551 | $g/mol$ |
| volumetric shape factor | $k_v$ | 1.112e-7 | $\ell/mm^3$ |
| saturation concentration | $c_s$ | 4.038 | $mol/\ell$ |
| crystal size distribution | $n(L,t)$ | | $\sharp/mm \cdot \ell$ |
| solute concentr. in liquid | $c(t)$ | | $mol/\ell$ |
| solute feed concentration | $c_f(t)$ | | $mol/\ell$ |

**Table 1.** Crystallizer parameters

Open-loop step responses of the nonlinear model are shown in Fig. 4.

Linearization about steady state with $n(L,t) = n_{ss}(L) + \Delta n(L,t)$, $c(t) = c_{ss} + \Delta c(t)$, $c_f(t) = c_{fss} + \Delta c_f(t)$, $\varepsilon(t) = \varepsilon_{ss} + \Delta\varepsilon(t)$, $\eta(t) = \eta_{ss} + \Delta\eta(t)$ leads to the linearized population balance

$$\Delta n_t = -k_g n'_{ss}(L)\Delta c - G(c_{ss})\Delta n_L - \frac{q}{V} h_{fp}(L)\Delta n \quad (20)$$

with initial condition $\Delta n(L,0) = 0$ and boundary condition

$$\Delta n(0,t) = \frac{3k_b}{k_g}(c_{ss} - c_s)^2 \Delta c(t), \quad (21)$$

and the linearized molar balance

$$\Delta c' = -\frac{q}{V}\Delta c + \frac{q}{V\varepsilon_{ss}}\Delta c_f + \frac{\rho - Mc_{ss}}{M\varepsilon_{ss}}\Delta\varepsilon' \quad (22)$$
$$+ \frac{q\rho - qMc_{fss} + q\rho z\eta_{ss}}{VM\varepsilon_{ss}^2}\Delta\varepsilon - \frac{q\rho z}{VM\varepsilon_{ss}}\Delta\eta$$

with $\Delta c(0) = 0$,

$$\Delta\varepsilon(t) = -k_v \int_0^\infty \Delta n(L,t)L^3 dL,$$

$$\Delta\eta(t) = k_v \int_{L_p}^\infty \Delta n(L,t)L^3 dL.$$

The infinite dimensional transfer function

$$G_{cry}(s) := \Delta c(s)/\Delta c_f(s)$$

is now computed formally as

$$G_{cry}(s) = \frac{p_{12}(s)}{p_{13}(s) + q_{12}(s)e^{-\frac{sL_f}{G(c_{ss})}} + r_{12}(s)e^{-\frac{sL_p}{G(c_{ss})}}}, \quad (23)$$

where $p_{12}, q_{12}, r_{12}, p_{13}$ are polynomials of order 12 respectively 13. In particular, $G_{cry}$ is meromorphic and strictly proper. If a class $\mathscr{K}$ of real rational proper controllers is used, hypotheses (i)–(iii) are satisfied, and (iv) is satisfied for $K$.

Before we apply our method, we verify hypothesis (iv) for $G$. We write the linearized system in the form

$$\begin{bmatrix} \Delta n_t \\ \Delta c' \end{bmatrix} = \begin{bmatrix} \mathscr{D} & \mathscr{M} \\ \mathscr{I} & \delta \end{bmatrix} \begin{bmatrix} \Delta n \\ \Delta c \end{bmatrix} + \begin{bmatrix} 0 \\ \gamma \end{bmatrix} \Delta c_f \quad (24)$$

where $\delta = -q/V - k_g \int_0^\infty n'_{ss}(L)L^3 dL$, and $\gamma = q/V\varepsilon_{ss}$ are constants, and $\mathscr{D} = -G(c_{ss})\frac{\partial}{\partial L} - (q/V)h_{fp}(L)$ is an unbounded differential operator on the Hilbert space

$$\mathscr{H} = L^2((0,\infty), \max\{1, L^3\}dL),$$

while $\mathscr{M} : \mathbb{R} \to \mathscr{H}$, $\Delta c \mapsto -k_g n'_{ss}(L)\Delta c$ is a multiplication operator, and $\mathscr{I} : \mathscr{H} \to \mathbb{R}$ is the bounded linear integral operator we obtain when we substitute the population balance equation to obtain $\Delta\varepsilon'(t) = -k_g\Delta c(t) \int_0^\infty n'_{ss}(L)L^3 dL + \int_0^\infty \frac{q}{V}h_{fp}(L)L^3\Delta n(L,t)dL + G(c_{ss})\int_0^\infty \Delta n_L(L,t)L^3 dL$. For the last term we use partial integration to obtain

$$-3G(c_{ss})\int_0^\infty \Delta n(L,t)L^2 dL,$$

so that altogether $\mathscr{I}[\Delta n] = \int_0^\infty \Delta n(L,t)\phi(L)dL$ for an expression $\phi(L)$ gathering weighted terms containing $L^2, L^3, h_{fp}(L)L^3$, and $\chi_{[L_p,\infty)}(L)L^3$ in (22). Setting $A = [\mathscr{D}, \mathscr{M}; \mathscr{I}, \delta]$, we have $D(A) = \{(u,v) \in \mathscr{H} \times \mathbb{R} : \frac{\partial u}{\partial L} \in \mathscr{H}, u(0) = (3k_b/k_g)(c_{ss} - c_s)^2 v\}$, we see that $A$ generates a strongly continuous semigroup of operators on a Hilbert space, while the input operator $B = [0; \gamma]$ is of finite rank. The same is true for the output operator $C = [0,1]$. It remains to check that $(A,B,C)$ is exponentially stabilizable and detectable.

**Lemma 5.** *The system* (20), (22) *with boundary condition* (21) *is exponentially stabilizable and detectable.*

*Proof.* For stabilizability the idea is to set up a linear integral operator $\mathscr{K}$ in state-feedback form $\frac{q}{V\varepsilon_{ss}}\Delta c_f(t) = \mathscr{K}[\Delta n, \Delta c](t)$ such that

$$\mathscr{K}[\Delta n, \Delta c] = -\frac{\rho - Mc_{ss}}{M\varepsilon_{ss}}\Delta\varepsilon'$$
$$- \frac{q\rho - qMc_{fss} + q\rho z\eta_{ss}}{VM\varepsilon_{ss}^2}\Delta\varepsilon + \frac{q\rho z}{VM\varepsilon_{ss}}\Delta\eta.$$

Because then substituting this control law in (22) leads to the equation $\Delta c' = -\frac{q}{V}\Delta c$, which is exponentially stable. Substituting $\Delta c$ back in (20) is then stable, because the differential operator $\mathscr{D} = -G(c_{ss})\frac{\partial}{\partial L} - \frac{q}{V}h_{fp}(L)$ with boundary condition (21) is exponentially stable. Checking boundedness of $\mathscr{K}$ is analogous to checking boundedness of the integral operator $\mathscr{I}$ above.

Concerning exponential detectability, in matrix notation the system may be written as

$$\begin{bmatrix} \Delta n_t \\ \Delta c' \end{bmatrix} = \left( \begin{bmatrix} \mathscr{D} & \mathscr{M} \\ \mathscr{I} & \delta \end{bmatrix} + \begin{bmatrix} \mathscr{G} \\ \eta \end{bmatrix} [0\ 1] \right) \begin{bmatrix} \Delta n \\ \Delta c \end{bmatrix}$$

where $\mathscr{I}, \mathscr{D}, \mathscr{M}, \delta$ are as in (24), $C = [0\ 1]$, and $\mathscr{F} = [\mathscr{G}; \eta]$ is sought. We choose $\mathscr{G} = -\mathscr{M}$, then the first equation becomes the exponentially stable $\Delta n_t = -G(c_{ss})\Delta n_L - \frac{4}{d}h_{fp}(L)\Delta n_L$, with boundary condition (21), which was already encountered in the previous proof. Substituting this back, the second equation becomes $\Delta c' = (\delta + \eta)\Delta c + r(t)$, which can be stabilized by choosing $\delta + \eta < 0$. That gives the required $\mathscr{F} = [\mathscr{G}; \eta]$. □

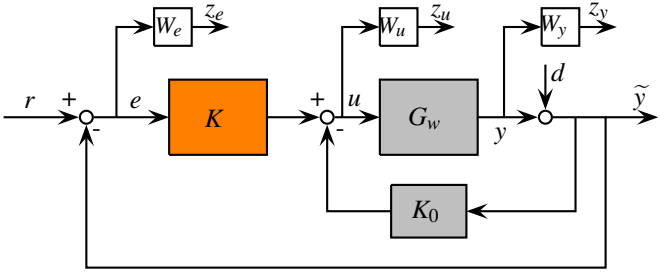

**Figure 5.** Control configuration for continuous crystallizer.

A first application of algorithm 1 reveals two unstable poles of $G_{\text{cry}}$. Using `systune` based on [1, 37, 58] we compute a static controller $K_0$ which stabilizes a low-order finite-difference crystallizer model $G_{502}$ with 502 states, where the target decay rate is chosen as 1e-7. Using algorithm 1, we then confirm that $K_0$ also stabilizes the infinite dimensional $G_{\text{cry}}(s)$.

In order to optimize performance of the continuous crystallizer, $G_{\text{cry}}$ is sampled as in algorithm 3, and the method is applied to $G = \mathscr{F}_l(G_{\text{cry}}, K_0)$, which has $n_p = 0$ rhp poles. We use the scheme of Fig. 5 with $K_0$ held fixed, while $K \in \mathscr{K}$ is optimized over the class $\mathscr{K}_{2,\text{stab}}$ of stable second-order controllers. The $H_\infty$-channel is $(r, d) \to (z_e, z_u, z_y)$ with weighing filters $W_e(s) = \frac{0.1s + 0.199}{s + 0.0199}$, $W_u = 0.01$, $W_y = \frac{100000s + 1.333e04}{s + 4.216e\text{-}06}$, and the optimal $H_\infty$-controller achieves a gain of $\gamma_\infty = 1.18$. The optimal controller of order 2 so obtained is

$$K^*(s) = \frac{54.47s^2 + 2.317s + 0.02446}{s^2 + 0.002033s + 4.374e\text{-}06},$$

and closed-loop stability is certified with algorithm 1.

As the last step the nonlinear crystallizer is simulated in closed loop with controller $K^*$. See Fig. 6.

The simulation uses a finite-difference semi-discretization with 4000 spatial steps. Blue shows controlled, red and magenta uncontrolled linearized and nonlinear state $c(t)$. The spatial resolution required for the desired precision is critical for a state-space control approach, and if state-space were used for control, system reduction would be inevitable.

# 8. Delay systems

Systems with delays are conveniently addressed by our novel synthesis technique. Semigroup theory is available [21], and the Nyquist test is applicable under hypothesis (iv). Standard tests for stabilizability and detectability exist and resemble those for rational systems. In the following, we illustrate the efficiency of our method in four typical studies.

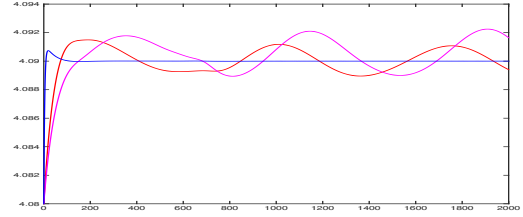

**Figure 6.** Simulation of $K^*$ with nonlinear system. The system is steered from old steady-state $c_{ss} = 4.08$ to new steady-state at $c_{ss} = 4.09$. Time is in minutes.

In the process industry, dead-time is a common phenomenon which may cause standard controllers to over-react to disturbances or set-point changes. The practical question is to decide whether or not dead-time is significant enough to be accounted for. One way to handle this is the celebrated Smith predictor [59] shown in Fig. 7. It applies to systems of the form $G(s) = G_0(s)e^{-\tau s}$, where $\tau$ is the delay, and where the delay-free $G_0(s)$ is called the lag. Typically, delay and lag are not precisely known, and we assume here for the purpose of illustration that a frequency sampled version $G(j\omega_v)$ of $G(s)$ is available for synthesis via (2). The Smith scenario now requires a model $G_m(s) = G_{m,0}(s)e^{-\tau' s}$ of the process, where $G_{m,0}$ is an estimation of the lag, $\tau'$ an estimation of the delay.

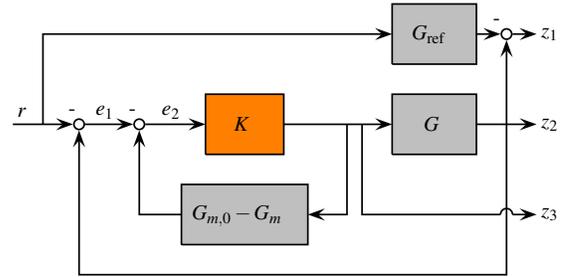

**Figure 7.** Synthesis interconnection with Smith predictor

**Example 4. Lag dominant process.** Our first delay study uses an example from [60], where it is assumed that the lag of $G$ is correctly identified, while an inaccurate guess $\tau' = 5.0$ of the true dead time $\tau = 5.5$ is made:

$$G(s) := G_{m,0}(s)e^{-5.5s}, \ G_m(s) := G_{m,0}(s)e^{-5.0s}. \quad (25)$$

We assume that frequency samples $G(j\omega_v)$ of the dead-time system $G$ are available on a grid $\Omega_{\text{opt}}$, but that on demand further sampled values $G(j\omega)$ can be obtained. Dead-time-free and delay-free reference systems are

$$G_{m,0}(s) := \frac{1}{(1 + 5s)(1 + 10s)}, \quad G_{\text{ref},0}(s) = \frac{3}{3 + 10s}.$$

Since the time constant of the process exceeds the delay time, the process is *lag-dominant*. The purpose of the reference model is to reduce the lag in closed loop, and it includes the incompressible model delay $\tau' = 5.0$, which leads to $G_{\text{ref}}(s) = G_{\text{ref},0}(s)e^{-\tau's}$.

The $H_\infty$-control problem minimizes the channel $r \to z = (z_1, z_2, z_3)$ with weight $W(s) = [w_1(s); w_2(s); w_3]$. Here $T_{z_1r}$ reflects set-point tracking, with $w_1(s) = \frac{0.01s+0.5986}{s+0.005986}$ a low-pass filter with crossover frequency at twice the bandwidth of the reference model $G_{\text{ref}}$. The channel $T_{z_2r}$ assesses mismatch between model $G_m$ and system $G$ on an appropriate frequency range described by the robustness weight $w_2(s) = \frac{2.39s+0.4078}{s+2.044}$, which is built as a tight upper bound of the relative uncertainty between the 3 models in (25). Finally, to limit the control effort the transfer function $T_{z_2r}$ with weight $w_3 = 0.1$ is included in the objective. Problem (2) is now solved via algorithm 3, where primary controllers $K$ are in the class $\mathcal{K}_{\text{pid}}$ of SISO PIDs. The optimal controller

$$K^*(s) = 2.93 + \frac{0.207}{s} + \frac{9.46s}{1+1.64s}, \qquad (26)$$

is obtained in 14s CPU using 33 iterations. It achieves good step responses, as seen in Fig. 8, and requires lower gain in the high frequency range compared to the loop-shaping controller given in [60], as seen in Fig. 9. $K^*$ is competitive with other controllers proposed for this study in the literature [61, 62, 63, 64].

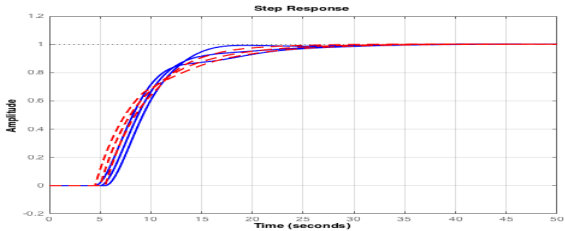

**Figure 8.** Step responses. Primary PID controller (26) (solid), controller in [60] (dotted)

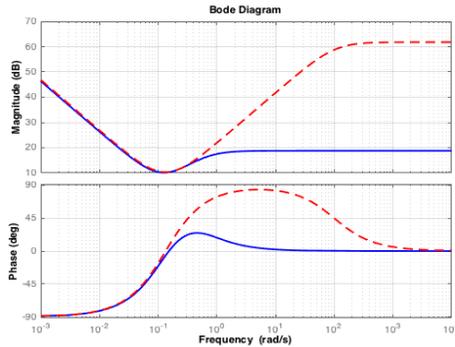

**Figure 9.** Bode plots. PID (solid), controller in [60] (dotted)

**Example 5. Dead-time dominant process.** Our second delay study is from [65] and follows again Fig. 7. With the same notations

$$G(s) = \frac{5.0}{1+38s}e^{-90s}, \qquad G_m(s) = \frac{5.6e^{-93.9s}}{1+40.2s},$$
$$G_{\text{test}}(s) = \frac{6}{1+42s}e^{-100s}, \qquad G_{\text{ref}}(s) = \frac{e^{-93.9s}}{1.33.33s}, \qquad (27)$$

which due to the large delay is now *dead-time dominant*. Here $G_{\text{test}}$ is used for posterior testing. Weighting filters are given as $w_3 = 0$ and

$$w_1(s) = \frac{5}{(20s+1)^2}, \, w_2(s) = \frac{2.661s+0.04519}{s+0.2265}.$$

The primary controller is a PI, and algorithm 3 gives the optimal $K_1 = 0.141 + 0.00645/s$. Step responses are shown in Figure 10 and exhibit significant overshoots and undershoots, which are chief features of long time-delay systems.

The transient behavior can be improved if larger settling times are accepted. With the modified reference model $G_{\text{ref}} = e^{-93.9s}/(1+70s)$ better transients are obtained, as seen in Fig. 10. The new primary PI obtained with algorithm 2 is now $K_2 = 0.0729 + 0.00322/s$.

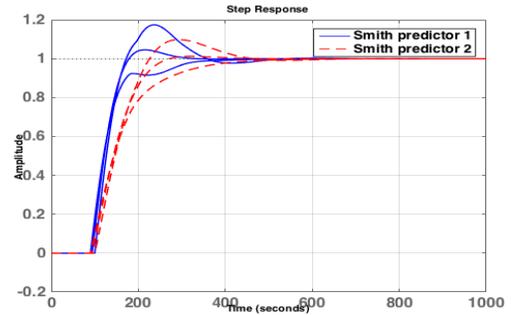

**Figure 10.** Step responses. PI primary controller $K_1$ (solid), PI primary controller $K_2$ (dashed)

**Example 6. Cavity flow.** A detailed study of cavity flows is given in [66, 67]. This challenging problem is taken from [68], where the infinite dimensional transfer function is available analytically as

$$G(s) = \frac{e^{-\tau_1 s}}{p_2(s) + q_2(s)e^{-\tau_2 s} + ce^{-\tau_3 s}},$$

with quadratic polynomials $p_2, q_2$. The $H_\infty$-objective is

$$\|(W_1 S, W_2 T)\|_\infty,$$

where $W_1(s) = (0.01s+502.5)/(s+50.25)$, $W_2(s) = (100s+500)/(s+50000)$. Optimization is over the class $\mathcal{K}_2$ of order 2 controllers. The optimal $K^*(s) = (0.718s^2+224.7s+2642)/(s^2+535.8s+2.268e04)$ achieves a final gain of $\gamma^* = 5.41$. The final grid size is $|\Omega_{\text{opt}}| = 382$, no update was necessary. Frequency responses are given in Fig. 11. Note that this test case can be approached using `systune` but requires a 15th-order Padé for the delay resulting in a 47th-order plant.

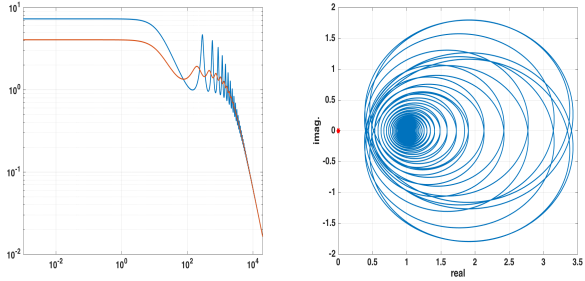

**Figure 11.** Cavity flow problem from [68]. Left image shows magnitude of $G(s)$ (blue) and $GS$ in closed loop (red). Right image shows the final Nyquist curve for $\omega \in [-3e4, \ 3e4]$.

**Example 7. Van de Vusse reactor [49].** An $H_\infty$ problem for a heat exchanger was solved in [49], and here this model is applied to a Van de Vusse reactor. Weighting filters were chosen as $W_u = 0.1$,

$$W_e = \frac{10^{-5}s + 1.502}{s + 0.07509}, \quad W_n = \frac{0.00125s^2 + 0.00035s + 5.10^{-5}}{2.5.10^{-5}s^2 + 0.007s + 1}$$

where $W_e$ and $W_u$ penalize tracking error and control effort, respectively. The filter $W_n$ specifies the frequency content of a noise input. While [49] considers full-order controllers of a suitable rational approximation, we use our transfer-function based approach from section 6, where we restrict for practical reasons optimization to the class $\mathscr{K}_3$ of 3rd-order controllers. Algorithm 3 yields an optimal $K^* \in \mathscr{K}_3$ with state-space representation

$$\left[ \begin{array}{ccc|c} -27.5666 & -26.2507 & 0 & 3.5532 \\ 21.9919 & -6.3124 & 2.9680 & 16.2390 \\ 0 & 0.6141 & -1.6018 & 2.2726 \\ \hline 4.4793 & -2.3704 & 1.8102 & -0.0768 \end{array} \right],$$

with certified locally optimal value $\gamma^* = 0.464$.

**Example 8. MIMO delay.** We consider studies from [69] with MIMO processes $G(s)$ with multiple input/output delays $G_{ij}(s) = G_{ij}^0(s)e^{-\tau_{ij}s}$, where $G_{ij}^0(s)$ is rational. All systems are square with dimensions 2 to 4.

The control scheme is a mixed sensitivity problem as in figure 5. For $G$ the first $2 \times 2$ study in [69] weightings are $W_e = w_e I_2$, $W_u = 0.01 I_2$, $W_y = w_y I_2$ and $w_e(s) = \frac{.01s + .2512}{s + .02512}$, $w_y(s) = \frac{100s + 5}{s + 500}$. The final $H_\infty$-norm is $\gamma^* = 1.07$, and is certified using Lemma 3. The method ends with $|\Omega_{\mathrm{opt}}| = 766$, for which it needs one update of the grid. The optimal $K^* \in \mathscr{K}_3$ was obtained as $[A_K \ B_K; \ C_K \ D_K] =$

$$\left[ \begin{array}{ccc|cc} -6.407 & 0.128 & 0 & -0.04964 & 3.273 \\ -1.189 & 0.006017 & -0.04487 & -0.5266 & 0.7091 \\ 0 & 0.02529 & -0.0794 & -0.439 & -0.2358 \\ \hline -2.454 & -0.1121 & 0.1447 & 0.2772 & 1.092 \\ -27.94 & 0.6764 & -0.4109 & -0.1795 & 15.78 \end{array} \right].$$

When allowed random restarts, `systune` with order 3 Padé approximation gives the same $H_\infty$-norm. Results for the re-maining $2 \times 2$, $3 \times 3$ and $4 \times 4$ examples from [69] are collected in Table 3 and the details are available upon request.

# 9. Comparison with convex-concave procedure

In this section we compare our approach to the convex-concave procedure (CCP) of [12, 13, 11]. The example is taken from [13], with process $G$ given as

$$G(s) = \begin{bmatrix} \frac{1}{s+1} & \frac{0.2}{s+3} & \frac{0.3}{s+0.5} \\ \frac{0.1}{s+2} & \frac{1}{s+1} & \frac{1}{s+1} \\ \frac{0.1}{s+0.5} & \frac{0.5}{s+2} & \frac{1}{s+1} \end{bmatrix}.$$

The problem is a standard mixed-sensitivity problem involving the weighted transfer functions $W_1S$ and $W_2KS$ with $W_1 := (s+3)/(3s+0.3)$ and $W_2 := (10s+2)/(s+40)$. The fine frequency grid $\Omega_{\mathrm{fine}}$ covers the interval $[10^{-2}, 10^2]$ with $N = |\Omega_{\mathrm{fine}}| = 1000$ points.

In [13], controllers are chosen as matrix fractions of polynomials $K(s) := N(s)D(s)^{-1}$, with

$$N(s) = N_d s^d + \ldots + N_1 s + N_0, \ D(s) = I_d s^d + \ldots + D_1 s + D_0.$$

Formally these $K$ have high order, but can be substantially reduced by taking minimal realization of order $\deg \det D(s)$. For instance, with $d = 3$ the fractional controller has order 27, but can be reduced to order 9. For comparison we compute controllers in state-space form (3) of increasing order $\mathrm{size}(A_K)$ using our approach and stop when no further progress is observed. Results are summarized in Table 2. To evaluate our non-smooth approach, we also report execution times (column 4) and the number of frequencies $|\Omega_{\mathrm{opt}}|$ that were used (column 5).

**Table 2.** Comparison of CCP procedure with trust-region non-smooth technique

| $K$ order | $\gamma^*$ CCP | $\gamma^*$ non-smooth | cpu time (sec.) | $|\Omega_{\mathrm{opt}}|$ |
|---|---|---|---|---|
| 1 | na* | 6.27 | 4.47 | 106 |
| 2 | na | 5.14 | 6.39 | 106 |
| 3 | 1.52 | 1.42 | 12.85 | 106 |
| 4 | na | 1.22 | 10.57 | 106 |
| 5 | na | 1.22 | 12.19 | 106 |
| 6 | 1.25 | 1.21 | 10.17 | 106 |
| 7 | na | | | |
| 8 | na | | | |
| 9 | 1.22 | | | |
| 15 | 1.21 | | | |

na: non available

Both techniques give comparable results. Our approach reaches the globally optimal value $\gamma^* = 1.21$ for a lower-order controller, taking advantage of working with state-space representations (3). Also, our algorithm uses a fairly small number of frequencies for both stability and performance, showing that sampling at higher densities is unnecessary. State-space data of the 4th-order optimal controller with certified $\gamma^* = 1.22$ is given as:

$$\begin{bmatrix} -1.1297 & 2.4687 & 0 & 0 & -0.62558 & 2.1069 & 2.4133 \\ 0.49873 & -1.7456 & -0.0056164 & 0 & -0.34424 & 1.25 & 0.07677 \\ 0 & 0.23265 & 0 & -0.17205 & -0.78512 & -2.9273 & 1.2731 \\ 0 & 0 & 0.047156 & -0.24021 & -0.86441 & -0.76551 & 0.29947 \\ \hline -0.2634 & -0.17576 & 0.41289 & -2.1632 & 0.11309 & -0.016873 & -0.01568 \\ -0.8761 & 2.5224 & -0.64724 & 0.59545 & 0.0138 & 0.11148 & -0.0338 \\ 1.0477 & -1.3105 & 0.43259 & -0.22984 & 0.0107 & 0.03145 & 0.098511 \end{bmatrix}$$

## 10. More exhaustive testing

Our method was tested on a bench of 40 plants, where algorithm 3 could be crosschecked. Table 3 shows examples from the Compleib collection [19], identified by their acronyms in column 1. As these examples (1-28) are finite dimensional, `systune` based on [1, 58] was used to compare with a standard structured $H_\infty$-synthesis [1]. For these tests the controller structure $\mathscr{K}_6$ of 6th-order controllers was used. Column 'algo. 3' gives the result of algorithm 3, with $|\Omega_{opt}|$ the size of the grid on exit, where column 'updates' gives the number of restarts in step 4 of algorithm 3. For instance, in study 'MFP' our method computed $K^* \in \mathscr{K}_6$ with optimal gain $\gamma^* = 4.27$ certified on exit. That is, $\|T_{wz}(P, K^*)\|_\infty = 4.27 + \vartheta$ with $|\vartheta| < 10^{-2}$. This was obtained with $|\Omega_{opt}| = 80$ and required 1 updating. Running `systune` on the same example gave $K_{sys}$ with the same structure and slightly better gain $\gamma_{sys} = 4.20$.

| test | systune | algo. 3 | $\Omega_{opt}$ | updates |
|------|---------|---------|------------|---------|
| AC3 | 3.10 | 2.98 | 2115 | 1 |
| AC6 | 3.52 | 3.65 | 108 | 1 |
| AC15 | 14.87 | 14.93 | 82 | 1 |
| AC16 | 14.86 | 14.87 | 90 | 1 |
| AC17 | 6.61 | 6.61 | 44 | 1 |
| HE2 | 2.45 | 2.45 | 112 | 2 |
| DIS1 | 4.16 | 4.17 | 105 | 1 |
| DIS3 | 1.04 | 1.05 | 185 | 1 |
| TG1 | 3.47 | 3.47 | 243 | 1 |
| AGS | 8.17 | 8.17 | 92 | 1 |
| WEC2 | 3.60 | 3.60 | 272 | 1 |
| WEC3 | 3.77 | 3.77 | 283 | 1 |
| BDT1 | 0.27 | 0.27 | 37 | 1 |
| MFP | 4.20 | 4.27 | 80 | 1 |
| UWV | 0.00 | 0.00 | 582 | 2 |
| EB1 | 3.09 | 3.10 | 140 | 2 |
| EB2 | 1.77 | 1.78 | 198 | 1 |
| PSM | 0.92 | 0.92 | 89 | 1 |
| NN4 | 1.29 | 1.29 | 118 | 1 |
| NN8 | 2.36 | 2.36 | 58 | 1 |
| NN11 | 0.0155 | 0.0255 | 138 | 1 |
| HF2D12 | 1037666.47 | 1037666.23 | 89 | 1 |
| HF2D13 | 101548.53 | 101548.54 | 211 | 1 |
| CM1 | 0.82 | 0.82 | 136 | 1 |
| DLR1 | 0.07 | 0.07 | 119 | 1 |
| JE3 | 4.14 | 4.15 | 958 | 1 |
| DLR2 | 201.28 | 147.22 | 6203 | 6 |
| DLR3 | 382.51 | 504.37 | 3546 | 7 |
| heat-N | – | 0.39 | 26 | 0 |
| heat-D | – | 0.60 | 17 | 1 |
| heat-M | – | 0.66 | 11 | 0 |
| reactor [49] | – | 0.46 | 101 | 1 |
| beam [70] | 0.14[70] | 0.14 | 201 | 1 |
| state-delay [71] | 0.2019 | 0.2015 | 81 | 0 |
| MIMO delay1 [69] | 1.07 | 1.07 | 766 | 1 |
| MIMO delay2 [69] | 1.61 | 1.61 | 260 | 1 |
| MIMO delay3 [69] | 1.59 | 1.48 | 195 | 0 |
| MIMO delay4 [69] | 0.47 | 0.51 | 2387 | 2 |
| cavity [68] | 5.55 | 5.41 | 382 | 0 |
| crystallizer | - | 1.18 | 550 | 0 |

**Table 3.** Test bench with 28 CompLeib examples and 12 infinite-dimensional studies

The 12 infinite-dimensional examples in Table 3 include in particular the studies heat-N, heat-D, heat-M, which use example 1 with Neumann, Dirichlet and mixed boundary conditions, where optimization is over the class $\mathscr{K}_1$ of first-order controllers. The $H_\infty$-controllers are $K_N(s) = (1.318s - 45.64)/(s + 4.493)$, $K_D(s) = (1.602s + 14.05)/(s + 0.2962)$ $K_M(s) = (5.885s + 12.31)/(s + 0.2916)$. In all heat studies the weights $W_e(s) = (0.01s + 3.015)/(s + 0.3015)$, $W_u = 0.01$ and $W_y(s) = (100s + 10)/(s + 1000)$ were used. The goal of each design was to track the set-point temperature at $\xi_0 = 1/3$, and to attenuate high frequency measurement noise.

The state-delay study uses a system with 2 states and state delay from [71]. The weights are $W_e(s) = (0.001s + 5.244)/(s + 0.5244)$, $W_u(s) = W_y(s) = (100s + 1.5)/(s + 1500)$, the channel is $r \to (W_e e, W_u u, W_y y)$. Here, $e$ denotes the tracking error with a reference model $\frac{3^2}{s^2 + 2 \times 0.8 \times 3s + 3^2}$. The optimal 2nd-order 2-DOF controller obtained by algorithm 3 is $K^*(s) = [0.9032s^2 + 7.546s + 9.488, -0.9052s^2 - 6.869s + 3.803]/(s^2 + 0.9293s + 6.63)$. As before, this result is certified with $\vartheta = 1e-2$ absolute accuracy, and crosschecked by `systune` with a 4th-order Padé approximation of the delay.

## 11. Conclusion

In this paper, we have presented a novel method for the synthesis of structured LTI controllers for a large class of infinite-dimensional systems described by their frequency response. Our method leverages non-smooth optimization techniques to compute locally optimal $H_\infty$-controllers.

Several frequency sampling techniques have been studied and a new adaptive sampling method for synthesis has been derived, which allows to certify exponential stability in closed loop and to computes $H_\infty$-performance of the resulting controllers within a fixed tolerance level $\vartheta$.

Our method is applicable to a fairly broad class of infinite-dimensional systems, including delay and integral-differential equations, boundary and distributed control of PDEs, and systems described by frequency-domain data. Local optimality certificates for program (2) are provided, and numerical testing confirms the excellent performance of the method, which often finds global optima. The method was evaluated on a large test bench including linearized PDEs, state-delayed and MIMO dead-time systems, and more detailed studies in process control.

## References

[1] P. Apkarian and D. Noll. Nonsmooth $H_\infty$ synthesis. *IEEE Trans. Automat. Control*, 51(1):71–86, January 2006.

[2] R. Hettich and K. O. Kortanek. Semi-infinite programming: theory, methods, and applications. *SIAM Review*, 35:380–429, 1993.

[3] S. Boyd and V. Balakrishnan. A regularity result for the singular values of a transfer matrix and a quadratically convergent algorithm for computing its $\mathbf{L}_\infty$-norm. *Syst. Control Letters*, 15:1–7, 1990.

[4] P. Apkarian, D. Noll, and L. Ravanbod. Nonsmooth bundle trust-region algorithm with applications to robust sta-

bility. *Set-Valued and Variational Analysis*, 24(1):115–148, 2016.

[5] P. Apkarian, D. Noll, and L. Ravanbod. Computing the structured distance to instability. In *SIAM Conference on Control and its Applications*, pages 423–430, 2015.

[6] E. Polak and Y. Wardi. A nondifferentiable optimization algorithm for the design of control systems subject to singular value inequalities over a frequency range. *Automatica*, 18(3):267–283, 1982.

[7] I. Horowitz. Quantitative feedback theory. *IEE Proc.*, 129-D(6):215–226, November 1982.

[8] E. van Solingen, J.W. van Wingerden, and T. Oomen. Frequency-domain optimization of fixed-structure controllers. *International Journal of Robust and Nonlinear Control*, 2016.

[9] Gorka Galdos, Alireza Karimi, and Roland Longchamp. $H_\infty$ controller design for spectral mimo models by convex optimization. *Journal of Process Control*, 20(10):1175–1182, 2010.

[10] A. Karimi, M. Kunze, and R. Longchamp. Robust PID controller design by linear programming. In *2006 American Control Conference*, pages 3931–3836, June 2006.

[11] M. Hast, K. J. Åström, B. Bernhardsson, and S. Boyd. PID design by convex-concave optimization. In *2013 European Control Conference (ECC)*, pages 4460–4465, July 2013.

[12] S. Boyd, M. Hast, and K.J. Åström. MIMO PID tuning via iterated LMI restriction. *International Journal of Robust and Nonlinear Control*, 26(8):1718–1731, 2016. rnc.3376.

[13] A Karimi and C Kammer. A data-driven approach to robust control of multivariable systems by convex optimization. *arXiv:1610.08776*, 2016.

[14] Thomas Lipp and Stephen Boyd. Variations and extension of the convex–concave procedure. *Optimization and Engineering*, 17(2):263–287, Jun 2016.

[15] H. Logemann. On the Nyquist criterion and robust stabilization for infinite-dimensional systems. In M. A. Kaashoek, J. H. van Schuppen, and A. C. M Ran, editors, *Robust Control of Linear Systems and Nonlinear Control: Proceedings of the International Symposium MTNS-89, Volume II*, pages 627–634. Birkhäuser Boston, Boston, MA, 1990.

[16] A. Sasane. An abstract Nyquist criterion containing old and new results. *Journal of Mathematical Analysis and Applications*, 370(2):703 – 715, 2010.

[17] M. Fardad and B. Bamieh. An extension of the argument principle and Nyquist criterion to a class of systems with unbounded generators. *IEEE Trans. Aut. Control*, 53(1):379–384, 2008.

[18] R. Curtain. A synthesis of time and frequency domain methods for the control of infinite-dimensional systems: A system theoretic approach. *SIAM Frontiers in Applied Mathematics*, 1989.

[19] F. Leibfritz. COMPL$_e$IB, COnstraint Matrix-optimization Problem LIbrary - a collection of test examples for nonlinear semidefinite programs, control system design and related problems. Technical report, Universität Trier, 2003.

[20] K. Zhou, J. C. Doyle, and K. Glover. *Robust and Optimal Control*. Prentice Hall, 1996.

[21] R. F. Curtain and H. Zwart. *An Introduction to Infinite-Dimensional Linear Systems Theory*, volume 21 of *Texts in Applied Mathematics*. Springer-Verlag, 1995.

[22] K.-J. Engel and R. Nagel. *One-Parameter Semigroups for Linear Evolution Equations*. Springer Graduate Texts in Math. Springer, 2000.

[23] S. G. Krantz. *Complex Variables: A Physical Approach with Applications and MATLAB*. Textbooks in Mathematics. Chapman and Hall/CRC, New York, 2007.

[24] D. Noll, O. Prot, and A. Rondepierre. A proximal bundle algorithm to minimize nonsmooth and nonconvex functions. *Pacific Journal of Optimization*, 4(3):569–602, 2008.

[25] D. Salamon. Infinite dimensional linear systems with unbounded control and observation: a functional analytic approach. *Transactions of the American Mathematical Society*, 300(2):383–431, 1987.

[26] G. Weiss. Transfer functions of regular systems. Part I: Characterizations of regularity. *Transactions of the American mathematical Society*, 342(2):827–854, 1994.

[27] R. Curtain. The Salamon-Weiss class of well-posed infinite-dimensional linear systems: a survey. *IMA Journal of Mathematical Control and Information*, 14:207 – 223, 1997.

[28] G. Weiss. Regular linear systems with feedback. *Mathematics of Control, Signals, and Systems*, 7(2):23–57, 1994.

[29] K.A. Morris. Justification of input-output methods for systems with unbounded control and observation. *IEEE-TAC*, 44(1):81–85, 1999.

[30] O.J. Staffans. *Well-Posed Linear Systems*. Encyclopedia of Mathematics and its Applications. Cambridge University Press, 2005.

[31] R. Rebarber. Conditions for the equivalence of internal and external stability for distributed systems with unbounded inputs. *IEEE Transactions on Automatic Control*, AC38:994–998, 1993.

[32] R. Curtain. Equivalence of input-output and exponential stability for infinite-dimensional systems. *Mathematica System Theory*, 21(4):1244–1265, 1988.

[33] S. Mossaheb. On the existence of right-coprime factorization for functions meromorphic in a half-plane. *IEEE Transactions on Automatic Control*, 25(3):550–551, Jun 1980.

[34] Hsiao-Ping Huang, Chung-Tarng Jiang, and Yung-Chen Chao. A new Nyquist test for the stability of control systems. *International Journal of Control*, 58(1):97–112, 1993.

[35] H. Zwart. Linearization and exponential stability. *arXiv:1404.3475v1*, 2014.

[36] S. Skogestad and I. Postlethwaite. *Multivariable Feedback Control - Analysis and Design*. Wiley, 1996.

[37] P. Apkarian and D. Noll. Nonsmooth optimization for multidisk $H_\infty$ synthesis. *European J. of Control*, 12(3):229–244, 2006.

[38] J.V. Burke, A.S. Lewis, and M.L. Overton. Two numerical methods for optimizing matrix stability. *Linear Algebra and its Applications 351-352*, pages 147–184, 2002.

[39] V. Bompart, P. Apkarian, and D. Noll. Nonsmooth techniques for stabilizing linear systems. In *Proc. American Control Conf.*, pages 1245–1250, New York, NY, July 2007.

[40] T. Kato. *Perturbation theory for linear operators; 2nd ed.* Grundlehren Math. Wiss. Springer, Berlin, 1976.

[41] A. Chang and K. Morris. Well-posedness of boundary control systems. *SIAM Journal of Control and Optimization*, 42(5):1101 – 1116, 2003.

[42] R. Curtain and G. Weiss. Well posedness of triples of operators (in the sense of linear system theory). In W. Schappacher F. Kappel, K. Kunisch, editor, *Control and Estimation of Distributed Parameter Systems*, pages 41–59. Birkhäuser Verlag, Basel, 1989.

[43] R. Curtain and K. Morris. Transfer functions of distributed parameter systems: A tutorial. *Automatica*, 45(5):1101 – 1116, 2009.

[44] G. Weiss and R. Rebarber. Dynamic stabilizability of well-posed linear systems. *5th International Symposium on Methods and Models in Automation and Robotics, Miedzyzdroje, Poland*, 1:2–9, 1998.

[45] K. Mikkola. State-feedback stabilization of well-posed linear systems. *Integral Equations and Operator Theory*, 55(2):249–271, 2006.

[46] K.A. Morris. $h_\infty$-output feedback of infinite-dimensional systems via approximation. *Systems and Control Letters*, 44(3):211 – 217, 2001.

[47] R. Triggiani. Boundary feedback stabilization of parabolic equations. *Appl. Math. Optim.*, 6(3):201–220, 1980.

[48] V. Barbu and R. Triggiani. Internal stabilization of the Navier-Stokes equations with finite-dimensional controllers. *Indiana University Mathematics Journal*, 53:1443–1494, 2004.

[49] H. Sano. $H_\infty$-control of a parallel-flow heat exchange process. *Bulletin of the Polish Academy of Sciences*, 65(1):11–19, 2017.

[50] W. Liu. Boundary feedback stabilization of an unstable heat equation. *SIAM J. Control Optim.*, 42(3):1033–1043, 2003.

[51] R. Curtain. The spectrum determined growth assumption for perturbations of analytic semigroups. *Systen and Control Letters*, 2(2):106 – 109, 1982.

[52] H. Logemann. On the transfer matrix of a neutral system: Characterizations of exponential stability in input-output terms. *System Control Letters*, 9:393–400, 1987.

[53] H. Zwart, Y. Le Gorrec, B. Maschke, and J. Villegas. Well-posedness and regularity of hyperbolic boundary control systems on a one-dimensional spatial domain. *ESAIM: Control, Optimisation and Calculus of Variations*, 16:1077–1093, 2010.

[54] M. Renardy. On the linear stability of hyperbolic PDEs and viscoelastic flows. *Zeitschrift für angewandte Mathematik und Physik ZAMP*, 45(6):854–865, 1994.

[55] J. Oostven. *Strongly stabilizable distributed parameter systems*. Frontiers in Applied Mathematics. SIAM, Philadelphia, 2000.

[56] U. Vollmer and J. Raisch. $H_\infty$-control of a continuous crystallizer. *Control Engineering Practice*, 9:837–845, 2001.

[57] A. Rachah, D. Noll, F. Espitalier, and F. Baillon. A mathematical model for continuous crystallization. *Mathematical Methods in the Applied Sciences*, 39(5):1101–1120, 2016.

[58] P. Apkarian, P. Gahinet, and C. Buhr. Multi-model, multi-objective tuning of fixed-structure controllers. In *European Control Conference (ECC)*, pages 856–861. IEEE, 2014.

[59] O. J. M. Smith. Closer control of loops with dead time. *Chemical Engineering Progress*, 53(9):217–219, 1957.

[60] Vinicius de Oliveira and Alireza Karimi. Robust Smith predictor design for time-delay systems with $H_\infty$ performance. *IFAC Proceedings Volumes*, 46(3):102 – 107, 2013.

[61] I. Kaya. Tuning Smith predictors using simple formulas derived from optimal responses. *Industrial & Engineering Chemistry Research*, 40(12):2654–2659, 2001.

[62] Z. J. Palmor and M. Blau. An auto-tuner for Smith dead time compensator. *Int. Journal of Control*, 60(1):117–135, 1994.

[63] T. Hagglund. A predictive PI controller for processes with long dead times. *IEEE Control Systems*, 12(1):57–60, Feb 1992.

[64] Chang-Chieh Hang, Qingguo Wang, and Li-Sheng Cao. Self-tuning Smith predictors for processes with long dead time. *International Journal of Adaptive Control and Signal Processing*, 9(3):255–270, 1995.

[65] P. Gahinet and L. F. Shampine. Software for modeling and analysis of linear systems with delays. In *Proc. American Control Conf.*, volume 6, pages 5600–5605, June 2004.

[66] P. Yan, M. Debiasi, X. Yuan, J. Little, H. Özbay, and M. Samimy. Experimental study of linear closed-loop control of subsonic cavity flow. *AIAA Journal*, 44(5):929–938, 2006.

[67] P. Yan, X. Yuan, H. Ozbay, M. Debiasi, E. Caraballo, M. Samimy, J. M. Myatt, and A. Serrani. Modeling and feedback control for subsonic cavity flows: A collaborative approach. In *Proceedings of the 44th IEEE Conference on Decision and Control*, pages 5492–5497, Dec 2005.

[68] Xin Yuan, Mehmet Önder Efe, and Hitay Özbay. On delay-based linear models and robust control of cavity flows. In Silviu-Iulian Niculescu and Keqin Gu, editors, *Advances in Time-Delay Systems*, pages 287–298. Springer Berlin Heidelberg, Berlin, Heidelberg, 2004.

[69] Qiang Xiong and Wen-Jian Cai. Effective transfer function method for decentralized control system design of multi-input multi-output processes. *Journal of Process Control*, 16(8):773 – 784, 2006.

[70] C. Foias, Hitay Özbay, and Allen R. Tannenbaum. *Robust Control of Infinite Dimensional Systems: Frequency Domain Methods*. Springer-Verlag New York, Inc., Secaucus, NJ, USA, 1996.

[71] M. Park, O. Kwon, J. Park, and S. Lee. Delay-dependent stability criteria for linear time-delay system of neutral type. *International Journal of Computer, Electrical, Automation, Control and Information Engineering*, 4(10):1602–1606, 2010.